\documentclass{amsart}

\newcommand{\figpage}[1]{\special{#1}\special{ps: /pageforpagesdotps exch def}}

\usepackage{graphics,verbatim, amsmath, amssymb, amsthm, amsfonts}	
\usepackage{ifthen}

\newtheorem{proposition}{Proposition}[section]
\newtheorem{theorem}[proposition]{Theorem}

\newtheorem{lemma}[proposition]{Lemma}

\newtheorem{cor}[proposition]{Corollary}

\theoremstyle{definition}

\theoremstyle{remark}

\numberwithin{equation}{section}

\usepackage[usenames]{color}

\newcommand{\margincolor}{Red}      

\newcounter{margincounter}
\newcommand{\marginnum}{\textcolor{\margincolor}{\begin{picture}(0,0)\put(5,3){\circle{13}}\end{picture}\arabic{margincounter}}}

\newcommand{\margin}[1]{\marginnum\marginpar{\textcolor{\margincolor}{\arabic{margincounter}.}\,\,\tiny #1}\addtocounter{margincounter}{1}}

\newcommand{\proofread}
{
\ifthenelse{\isundefined{\margin}}
{
\special{!userdict begin /bop-hook{1.2 1.2 scale -51 -60 translate}def end}
}
{
\special{!userdict begin /bop-hook{1.2 1.2 scale -91 -60 translate}def end}
\setboolean{@mparswitch}{false} 
\addtolength{\marginparwidth}{-3mm}
}
}
\proofread

\newcommand{\reals}{\mathbb R}

\newcommand{\set}[1]{{\lbrace #1 \rbrace}}

\renewcommand{\th}{^\mathrm{th}}

\newcommand{\join}{\vee}

\newcommand{\meet}{\wedge}

\newcommand{\covered}{\lessdot}
\newcommand{\A}{\mathcal{A}}
\newcommand{\F}{\mathcal{F}}
\newcommand{\E}{\mathcal{E}}
\newcommand{\pidown}{\pi_\downarrow}

\newcommand{\piup}{\pi^\uparrow}

\newcommand{\br}[1]{\langle #1 \rangle}

\newlength{\lsmash} 
\newlength{\myheight}
\newlength{\mydepth}

\newcommand{\up}[1]{
\settoheight{\myheight}{\ensuremath{\overline{#1}}}
\addtolength{\myheight}{-\lsmash}
\overline{
\protect \raisebox{0 pt}[\myheight][0 pt]{\ensuremath{\overline{#1}}}
}
}

\newcommand{\down}[1]{
\settodepth{\mydepth}{\ensuremath{\underline{#1}}}
\addtolength{\mydepth}{-\lsmash}
\underline{
\protect \raisebox{0 pt}[0 pt][\mydepth]{\ensuremath{\underline{#1}}}
}
}

\author{Nathan Reading}
\address{Department of Mathematics, North Carolina State University, Raleigh, NC, USA}


\title{From the Tamari lattice to Cambrian lattices and beyond}

\begin{document}

\begin{abstract}
In this chapter, we trace the path from the Tamari lattice, via lattice congruences of the weak order, to the definition of Cambrian lattices in the context of finite Coxeter groups, and onward to the construction of Cambrian fans.
We then present sortable elements, the key combinatorial tool for studying Cambrian lattices and fans.
The chapter concludes with a brief description of the applications of Cambrian lattices and sortable elements to Coxeter-Catalan combinatorics and to cluster algebras.
\end{abstract}

\maketitle


\section{A map from permutations to triangulations}\label{permtri sec}
The road from the Tamari lattice to Cambrian lattices starts with a simple map from the set $S_{n+1}$ of permutations of $\set{1,\ldots,n+1}$ to the set of triangulations of a convex polygon with $n+3$ vertices.
This map connects the Tamari lattice to the weak order on permutations, and opens the door to understanding the Tamari lattice in a broader lattice-theoretic context.

\begin{figure}
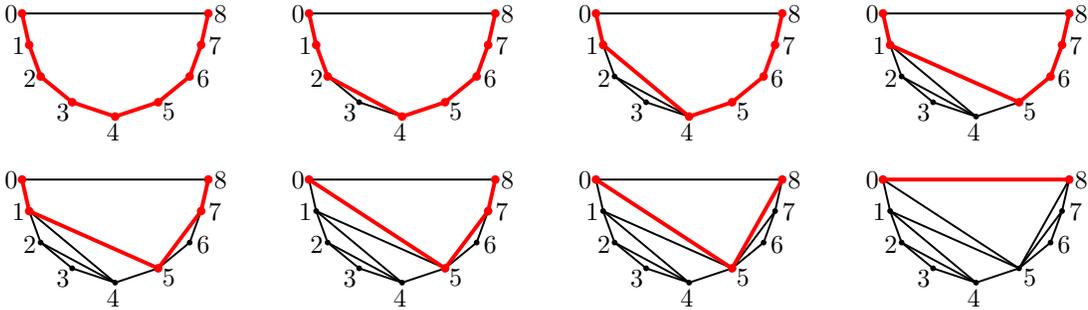

\begin{tabular}{ccccccc}
\figpage{ps: 1}\scalebox{.45}{\includegraphics{map.ps}}&&
\figpage{ps: 2}\scalebox{.45}{\includegraphics{map.ps}}&&
\figpage{ps: 3}\scalebox{.45}{\includegraphics{map.ps}}&&
\figpage{ps: 4}\scalebox{.45}{\includegraphics{map.ps}}\\[5 pt]
\figpage{ps: 5}\scalebox{.45}{\includegraphics{map.ps}}&&
\figpage{ps: 6}\scalebox{.45}{\includegraphics{map.ps}}&&
\figpage{ps: 7}\scalebox{.45}{\includegraphics{map.ps}}&&
\figpage{ps: 8}\scalebox{.45}{\includegraphics{map.ps}}
\end{tabular}
\caption{The triangulation $\eta(3246175)$}
\label{map fig}
\end{figure}

One of the many realizations of the Tamari lattice is as a partial order on triangulations of a convex polygon.  
Specifically, take $Q$ to be a convex $(n+3)$-gon in the plane and identify the vertices of $Q$ with the numbers $0,1,\ldots,n+1,n+2$.
We require that the vertices $0$ and $n+2$ be on a horizontal line, with $0$ to the left and with all other vertices below that line.
Furthermore, we require that the vertices $1$ through $n+1$ be placed so that, for all $i$ from $0$ to $n+1$, the vertex $i$ is strictly further left than the vertex $i+1$.
A correct construction of $Q$, for $n=7$, is shown in the top-left picture of Figure~\ref{map fig} for the case $n=7$.

\begin{figure}
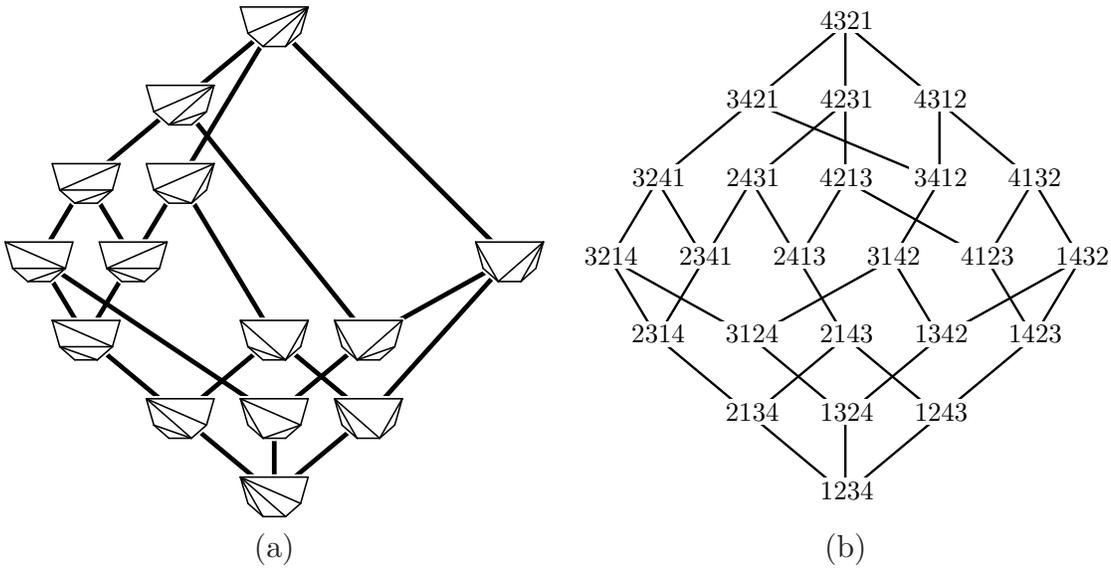

\begin{tabular}{cc}
\scalebox{.985}{\includegraphics{permtribottoms.ps}}
&
\scalebox{.985}{\includegraphics{weakS4.ps}}\\
(a)&(b)
\end{tabular}
\caption{a: The Tamari lattice.  b: The weak order on permutations.}
\label{tam-weakS4 fig}
\end{figure}

A \emph{triangulation} of $Q$ is a tiling of $Q$ by triangles whose vertices are contained in the vertex set of $Q$.
The triangulation is specified by the collection of $n$ diagonals of $Q$ appearing as edges of the triangles.
A \emph{diagonal flip} on a triangulation of $Q$ is the operation of removing a diagonal of the triangulation to create a quadrilateral from two triangles, and then inserting the other diagonal of the quadrilateral to create a new triangulation.
The Tamari lattice is a partial order on triangulations of $Q$ whose cover relations are given by diagonal flips.
The two triangulations in the cover differ by exactly one diagonal of $Q$, and the higher triangulation in the cover relation is the one in which this diagonal has larger slope.
The Tamari lattice, for $n=3$, is shown in Figure~\ref{tam-weakS4 fig}.

This definition of the Tamari lattice highlights its connection to the associahedron.
Since the vertices of the associahedron can be labeled by triangulations of a fixed convex polygon such that edges are given by diagonal flips, the Hasse diagram of the Tamari lattice is isomorphic to the $1$-skeleton of the associahedron.

To define the weak order, we first write permutations in one-line notation, meaning that we represent a permutation $x$ of $\set{1,\ldots,n+1}$ by the sequence $x_1x_2\cdots x_{n+1}$, where $x_i$ means $x(i)$.
There is a cover relation $x\covered y$ in the weak order whenever the one-line notations of $x$ and $y$ differ only by swapping a pair of adjacent entries.
The permutation $x$ is the one in which the two entries appear in numerical order, and $y$ is the permutation in which the two entries appear out of order.
For example, the weak order on $S_4$ is shown in Figure~\ref{tam-weakS4 fig}.b.

\begin{figure}
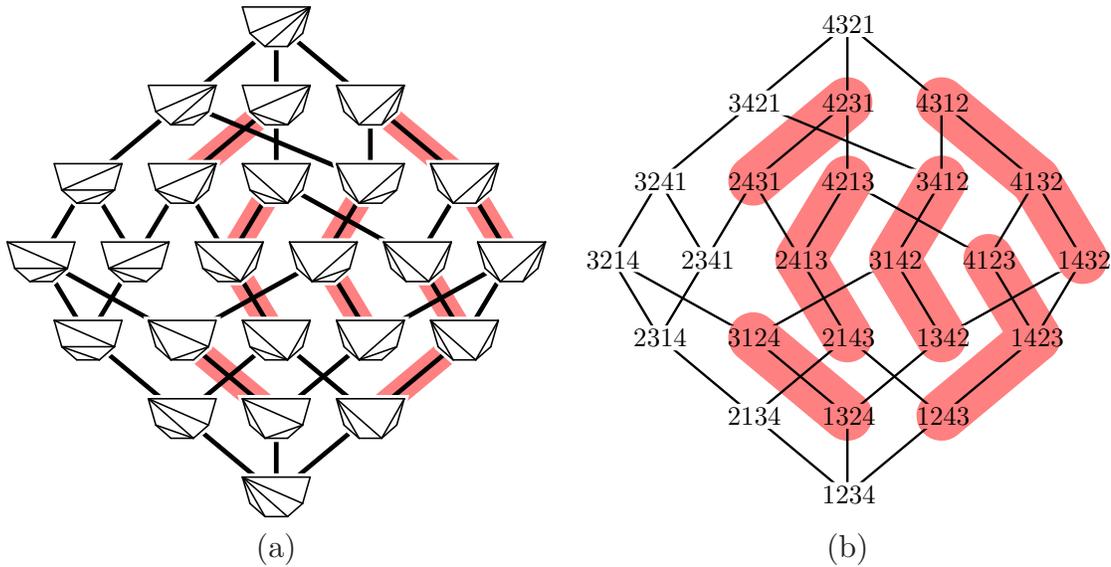

\begin{tabular}{cc}
\scalebox{.985}{\includegraphics{permtri.ps}}
&
\scalebox{.985}{\includegraphics{tamcong.ps}}\\
(a)&(b)
\end{tabular}
\caption{a:  The map $\eta$ applied to every permutation in $S_4$.  b: The Tamari congruence on $S_4$.}
\label{permtri-tamcong fig}
\end{figure}

We now define a map $\eta$ from $S_{n+1}$ to the set of triangulations of $Q$.
Start with a path along the bottom edges of $Q$, as shown in the first frame of Figure~\ref{map fig}.
Given a permutation $x\in S_{n+1}$, read from left to right in the one-line notation for $x$. For each entry, create a new path by deleting the corresponding vertex from the old path.
The triangulation $\eta(x)$ is defined by the union of the sequence of paths, as illustrated in Figure~\ref{map fig} for the permutation with one-line notation $3246175$.
Figure~\ref{permtri-tamcong fig}.a shows the result of applying $\eta$ to every permutation in $S_4$.
The shaded edges indicate covering pairs in the weak order which map to the same triangulation.

This map and similar maps have appeared in many papers, including \cite{Iterated,NonpureII,LR,LRorder,cambrian,Tonks}.   
The map can be seen in a broader context in the chapter by Rambau and Reiner \cite{RamRei} in this volume, specifically by giving some thought to \cite[Theorem~9]{RamRei} and the accompanying figure.

Bj\"{o}rner and Wachs \cite[Section~9]{NonpureII} studied a map $\tau$ from permutations to binary trees that is, up to a standard bijection from triangulations to binary trees, identical to $\eta$.
We describe their results in terms of the map $\eta$.
First, the fiber $\eta^{-1}(\Delta)$ of each triangulation $\Delta$ is a non-empty interval in the weak order on $S_{n+1}$.
A permutation is the minimal element in its $\eta$-fiber if and only if it avoids the pattern $312$.
That is, a permutation $x$ is minimal in its fiber if and only if there is no sequence of three (not-necessarily adjacent) entries in the one-line notation for $x$ such that the largest of the three is first, followed by the smallest of the three, and finally the median-valued.
Similarly, a permutation is the maximal element in its $\eta$-fiber if and only if it avoids the pattern $132$.
For example, comparing Figures~\ref{tam-weakS4 fig}.b and~\ref{permtri-tamcong fig}.a, we see that the permutation $4213$ is not the minimal element of its $\eta$-fiber, and indeed, the sequence $413$ (or the sequence $423$) is an instance of the pattern $312$ in the permutation $4213$.
However, $4213$ is the maximal element in its $\eta$-fiber because it avoids the pattern $132$.

Bj\"{o}rner and Wachs also showed that the weak order and the Tamari lattice are closely related.
Specifically, the restriction of the weak order to $312$-avoiding permutations is a sublattice of the weak order, and the restriction of $\eta$ to this sublattice is an isomorphism from the sublattice to the Tamari lattice.
This is readily seen in the case of $S_4$ by inspection of Figures~\ref{tam-weakS4 fig} and~\ref{permtri-tamcong fig}.a.
The sublattice of the weak order consisting of $312$-avoiding permutations (and thus the Tamari lattice) is also a quotient of the weak order in an order-theoretic sense.
Indeed, the results of \cite{NonpureII} go most of the way to establishing something stronger:
As we will see in Section~\ref{cong sec}, the Tamari lattice is a lattice quotient (i.e.\ a lattice-homomorphic image) of the weak order, because the map $\eta$ is a lattice homomorphism.
This is the key insight that leads to the notion of Cambrian lattices.

Before we shift the discussion to lattice theory, we give a generalization, in a more combinatorial direction, of the map $\eta$.
We will see in Section~\ref{camb sec} that this generalization is also an essential step towards Cambrian lattices.
The generalization, which was exploited in \cite{cambrian}, draws on the description in \cite[Section~4.3]{Equivariant} of a similar family of maps in the context of signed permutations.
These families of maps arise quite naturally in the context of (equivariant) iterated fiber polytopes, as explained in \cite{Iterated} and \cite[Section~4.3]{Equivariant} and as summarized in \cite[Sections~4, 6]{cambrian}.

\begin{figure}
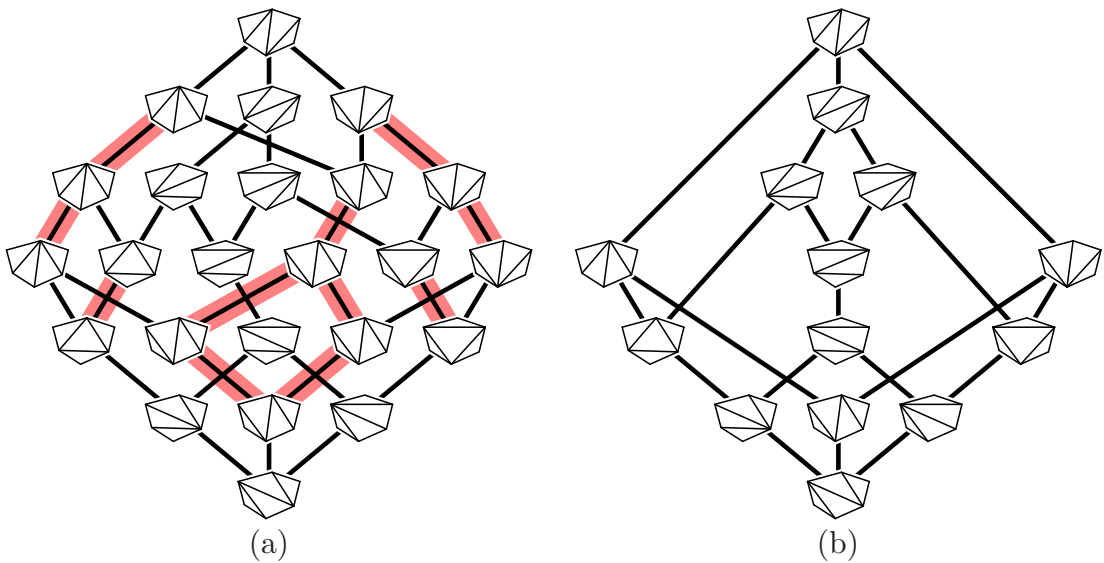

\begin{tabular}{cc}
\scalebox{.97}{\includegraphics{permtri2.ps}}
&
\scalebox{.97}{\includegraphics{permtri2bottoms.ps}}\\
(a)&(b)
\end{tabular}
\caption{a: A non-Tamari permutations-to-triangulations map applied to every permutation in $S_4$.  b:  A non-Tamari Cambrian lattice} 
\label{permtri2-nontam fig}
\end{figure}

To generalize $\eta$, we alter the construction of the polygon $Q$ by removing the requirement that the vertices $1$ through $n+1$ be located below the horizontal line containing $0$ and $n+2$.
We keep the requirement that, for all $i$ from $0$ to $n+1$, the vertex $i$ is strictly further left than the vertex $i+1$.
Again we start with a path along the bottom edges of $Q$, and read the one-line notation of a permutation from left to right.
When we read an entry whose corresponding vertex is on the bottom of $Q$, we \textbf{remove} that vertex from the path, as before.
When we read an entry whose corresponding vertex is on the top of $Q$, we \textbf{insert} that vertex into the path.
Figure~\ref{permtri2-nontam fig}.a shows this new permutations-to-triangulations map applied to all of the permutations in $S_4$, in the case where the vertices $1$, $2$, and $4$ are on the bottom of $Q$ and $3$ is on the top.
To avoid a profusion of notation, we use the symbol $\eta$ to refer to any of the permutations-to-triangulations maps, tacitly assuming a choice of $Q$.
We use the phrase ``the Tamari case'' to distinguish the original definition of $Q$ and of $\eta$.

As another example, consider the case where all of the vertices $1$ through $n$ are \textbf{above} the line containing $0$ and $n+2$.
In this case, the symmetries of the problem imply that $\eta$ has the same pattern-avoidance properties as described above, except that ``$312$'' is replaced by ``$231$'' throughout the description, and ``$132$'' is replaced by ``$213$'' throughout.

When some vertices are on top of $Q$ and others are on bottom, as in the example of Figure~\ref{permtri2-nontam fig}.a, the behavior of the map is a mixture of the ``$231$-behavior'' and the ``$312$-behavior,'' as we now explain.
The locations, top or bottom, of the vertices are recorded by upper- or lower-barring the symbols from $1$ to $n+1$.
Thus, for example, we write $\up{3}$ to indicate that the vertex $3$ is on top of the polygon $Q$ or we write $\down{3}$ to indicate that $3$ is on the bottom of $Q$.
In \cite[Proposition~5.7]{cambrian}, it is shown that a permutation is a minimal element in its $\eta$-fiber if and only if it avoids the patterns $31\down{2}$ and $\up{2}31$.
That is, the permutation may contain no $312$-pattern such that the ``$2$'' in the pattern is a lower-barred number and may contain no $231$-pattern such that the ``$2$'' in the pattern is an upper-barred number.

The subposet consisting of permutations avoiding $31\down{2}$ and $\up{2}31$ is a sublattice, isomorphic to a ``Tamari-like'' lattice on triangulations \cite[Theorem~6.5]{cambrian}.  
As in the original Tamari lattice, cover relations are diagonal flips and moving up in the order means increasing the slope.
In particular, just as in the Tamari lattice, the undirected Hasse diagram of the lattice is isomorphic to the $1$-skeleton of the associahedron.

Figure~\ref{permtri2-nontam fig}.b shows this Tamari-like lattice in the case where the vertices $1$, $2$, and $4$ are on the bottom of $Q$ and $3$ is on the top.  
These Tamari-like lattices are the first examples of Cambrian lattices beyond the Tamari lattices.
The next several sections develop Cambrian lattices in general.

\section{Bringing lattice congruences into the picture}\label{cong sec}
A lattice is a partially ordered set such that any pair $x,y$ of elements has a unique maximal lower bound (the meet of $x$ and $y$, written $x\meet y$) and a unique minimal upper bound (the join of $x$ and $y$, written $x\join y$). 
The Tamari lattice is, as its name suggests, a lattice, and so is the weak order.

The meet and join operation in a lattice satisfy certain properties, and indeed the definition of a lattice can be rephrased in universal-algebraic language: 
A lattice is a set with two operations $\meet$ and $\join$ satisfying a certain list of axioms.
The exact form of the axioms is not important for our present purposes, but the idea that a lattice is an algebraic object is fundamental.
From the algebraic viewpoint, it becomes natural to consider \emph{lattice homomorphisms} and \emph{lattice congruences}.
In this section, we point out that the permutations-to-triangulations map $\eta$ is a lattice homomorphism and characterize its associated congruence in a way that can be generalized to all finite Coxeter groups.

A lattice homomorphism is a map $\eta$ between lattices such that $\eta(x\meet y)=\eta(x)\meet\eta(y)$ and $\eta(x\join y)=\eta(x)\join\eta(y)$.
A lattice congruence is an equivalence relation $\Theta$ on a lattice $L$ such that, if $x_1\equiv y_1$ and $x_2\equiv y_2$ modulo $\Theta$, then $(x_1\meet x_2)\equiv(y_1\meet y_2)$ and $(x_1\join x_2)\equiv(y_1\join y_2)$ modulo $\Theta$.
Given a lattice congruence $\Theta$ on $L$, we construct the quotient lattice $L/\Theta$, whose elements are the $\Theta$-classes.
If $x$ and $y$ are elements of the lattice $L$, then the meet of the $\Theta$-class of $x$ and the $\Theta$-class of $y$ is the $\Theta$-class of $x\meet y$.
The join is defined similarly, and these operations are well-defined precisely because $\Theta$ is a congruence.
The nonempty fibers of a lattice homomorphism are the classes of a lattice congruence, and, given a lattice congruence on $L$, there is a natural lattice homomorphism from $L$ to the quotient lattice.

Some readers may only be familiar with the notion of congruences and quotients in the context of groups or rings, and perhaps without reference to the term ``congruence.''
In a group, the congruence classes are the cosets of a normal subgroup;  in a ring, they are the additive cosets of an ideal.
In a lattice, congruences are more complicated, but are still often amenable to study.

There is much to be gained, for the study of lattice congruences as for other lattice-theoretic topics, by passing back and forth between a universal-algebraic and an order-theoretic point of view.
This is especially true for finite lattices.
When $L$ is a finite lattice, congruences in $L$ have a particularly simple order-theoretic characterization.
It is straightforward to verify that an equivalence relation $\Theta$ on a finite lattice $L$ is a lattice congruence if and only if it satisfies the following three properties:
\begin{enumerate}
\item[(i) ] Every equivalence class is an interval.
\item[(ii) ] The projection $\pidown:L\rightarrow L$, mapping each element $a$ of $L$ to the minimal element in its $\Theta$-class, is order-preserving.
\item[(iii) ] The projection $\piup:L\rightarrow L$, mapping each element $a$ of $L$ to the maximal element in its $\Theta$-class, is order-preserving.
\end{enumerate}
This simple observation and its generalizations appear to be the subject of many independent rediscoveries, including in~\cite{Cha-Sn,Dorfer,Kolibiar,dissective}.

The order-theoretic rephrasing of the definition of a lattice congruence shows, in particular, how to recognize, combinatorially, when a map of lattices might be a lattice homomorphism.
For example, what we have already learned about the map $\eta$, in the Tamari case, leads us to suspect that the fibers of $\eta$ constitute a lattice congruence of the weak order, and thus to suspect that $\eta$ is a lattice homomorphism from the weak order to the Tamari lattice.
Indeed, our suspicions are correct, as is easily proven using the tools already developed in \cite{NonpureII};  the fibers of $\eta$, in the Tamari case, define a lattice congruence called the \emph{Tamari congruence}.
The Tamari congruence on $S_4$ is shown in Figure~\ref{permtri-tamcong fig}.b.
As a special case of the Cambrian construction, we will see that $\eta$ is a lattice homomorphism for the other choices of $Q$ as well.

The permutations-to-triangulations setup does more than simply find a lattice homomorphism between two given lattices.
Suppose we are given a surjective map $\eta$ from a lattice $L$ to a \textbf{set} $U$, and suppose that the fibers of $\eta$ satisfy conditions (i)--(iii).
The fibers therefore constitute a lattice congruence $\Theta$ on $L$, so there is a quotient lattice $L/\Theta$ such that the restriction of $\eta$ is a bijection from $L/\Theta$ to $U$.
We can use this bijection to \textbf{define} a lattice structure on $U$ isomorphic to $L/\Theta$, and when we do so, the map $\eta$ is a lattice homomorphism from $L$ to $U$.
In particular, the map $\eta$ from the weak order on permutations to the \textbf{set} of triangulations can be thought of as \textbf{defining} the Tamari lattice.
In \cite[Sections~4--5]{cambrian}, the Tamari lattice is constructed from the ground up by this method, together with the other Tamari-like lattices.

Indeed, to define the Tamari lattice, we don't even need the map $\eta$.
We can define the Tamari lattice to be the quotient of the weak order modulo the Tamari congruence.
In doing so, we lose the combinatorial realization of the Tamari lattice as a partial order on the set of triangulations.
But we gain a point of view on the Tamari lattice that leads us to a broad generalization.
To find this generalization, we need to understand what makes the Tamari congruence special among the vast number of lattice congruences of the weak order.

First, we describe some general features of congruences of a finite lattice $L$.
Much of this description does not generalize to infinite lattices, so we keep strictly to the finite case.
Suppose $\Theta$ is a congruence on $L$.  
If $x\covered y$ (i.e.\ $x$ is covered by $y$) in $L$ and $x\equiv y$ modulo $L$, then we say $\Theta$ \emph{contracts} the edge $x\covered y$.
Since congruence classes are intervals, we can describe $\Theta$ completely just by recording which edges in the Hasse diagram are contracted by $\Theta$.
Indeed, this is exactly how congruences are indicated in Figures~\ref{permtri-tamcong fig} and~\ref{permtri2-nontam fig}.a.

\begin{figure}
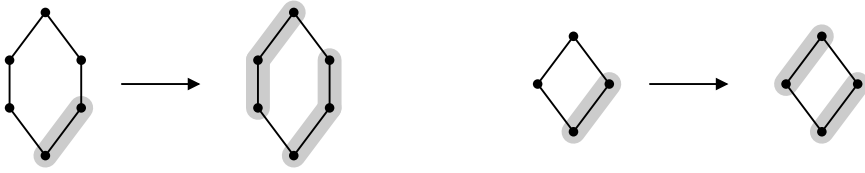

\scalebox{.3}{\includegraphics{move1.ps}}\qquad\qquad
\scalebox{.3}{\includegraphics{move2.ps}}
\caption{Local forcing requirements}
\label{move fig}
\end{figure}

As one would expect, edges cannot be contracted independently to form a congruence;  rather, contracting one edge may force other edges to be contracted.
As a simple example of this \emph{edge-forcing}, consider the hexagonal poset shown at the left of Figure~\ref{move fig}.
We call the two edges incident to the minimal element the \emph{bottom edges}, the two edges incident to the top element the \emph{top edges}, and the other two edges \emph{side edges}.
The following facts are easily checked using either the definition or the order-theoretic characterization of a lattice congruence:
For either side edge, there is a congruence contracting only that edge.
For any bottom or top edge, a congruence contracting that edge must also contract the opposite edge and both side edges.
Here the opposite edge means the edge related by a half-turn of the diagram.

In a general finite lattice, edge-forcing can be much more complicated.
In the weak order on $S_{n+1}$, however, edge forcing can be understood in terms of local forcing requirements of the kind illustrated in the hexagon example.
These intervals, and the local forcing requirements, are illustrated in Figure~\ref{move fig}.
In these pictures, shaded edges represent edges in the lattice that are contracted by a given congruence.  
The pictures show, given a contracted edge, which other edges must be contracted.
By symmetry (turning intervals upside-down or reflecting them in a vertical line), each requirement pictured describes four requirements.
In words, the local forcing requirements are that whenever a bottom or top edge of the interval is contracted, the opposite bottom or top edge is also contracted, as well as all side edges (if side edges are present).
A collection of edges is \emph{closed under local forcing} if the local forcing requirements are satisfied on every hexagonal and quadrilateral interval.
As a consequence of \cite[Theorem~25]{hyperplane}, we have the following characterization of lattice congruences of the weak order on permutations:
Given a collection $\E$ of edges in the weak order on $S_{n+1}$, there exists a lattice congruence $\Theta$ contracting exactly the edges in $\E$ if and only if $\E$ is closed under local forcing.

Using the local forcing requirements, the reader can readily verify that the Tamari congruence on $S_4$, shown in Figure~\ref{permtri-tamcong fig}.b, is the unique finest (in the sense of refinement of set partitions) lattice congruence contracting the edges $1324\covered3124$ and $1243\covered1423$.
This characterization generalizes to larger values of $n$:
In general, the Tamari congruence is the finest congruence contracting the edges 
\begin{multline*}1\cdots (j-2)(j-1)(j+1)j(j+2)\cdots (n+1)\\
\covered\,1\cdots (j-2)(j+1)(j-1)j(j+2)\cdots (n+1)
\end{multline*}
for $j=2,\ldots,n$.
If we vary the polygon $Q$, to alter the map $\eta$, the corresponding congruence $\Theta$ retains a similar characterization \cite[Theorem~6.2]{cambrian}.
It is the finest congruence contracting the edges 
\begin{multline*}
1\cdots (j-2)(j-1)(j+1)j(j+2)\cdots (n+1)\\
\covered\,1\cdots (j-2)(j+1)(j-1)j(j+2)\cdots (n+1)\end{multline*}
for $j=2,\ldots,n$ such that $j$ is lower-barred and the edges 
\begin{multline*}1\cdots (j-2)j(j-1)(j+1)(j+2)\cdots (n+1)\\
\covered\,1\cdots (j-2)j(j+1)(j-1)(j+2)\cdots (n+1)
\end{multline*} 
for $j=2,\ldots,n$ such that $j$ is upper-barred.

Although it is the above characterization that leads to the Cambrian lattices, it is worth mentioning another characterization of the Tamari congruence and its cousins, found recently by Santocanale and Wehrung~\cite{SW}.
The set of all congruences on a lattice is itself a finite (distributive) lattice under the refinement partial order.
The Tamari-like congruences are exactly the minimal \emph{meet-irreducible congruences} 
in the lattice of congruences of the weak order on permutations.
The natural generalization (into the context of finite Coxeter groups) suggested by the Santocanale-Wehrung characterization leads, not to the Cambrian lattices, but to a different object that has not yet been studied.

\section{Cambrian lattices}\label{camb sec}
The symmetric group $S_{n+1}$ is a classical example of a finite \emph{Coxeter group}.
Coxeter groups are given by a simple combinatorial presentation, which we review below, but there is also a much more geometric point of view:  A finite group can be realized as a Coxeter group if and only if it can be realized in Euclidean space as a group of transformations generated by reflections.

\emph{Generalized associahedra} were introduced, in combinatorial terms, by Fomin and Zelevinsky~\cite{ga} and were realized as polytopes by those authors and Chapoton~\cite{gaPoly}.
There is a generalized associahedron for each finite Coxeter group $W$, called the $W$-associahedron.
More precisely, there is a generalized associahedron for each \emph{root system}, but we gloss over this distinction for the purpose of brevity.
In fact, we skip the definition of the generalized associahedron entirely.
The key point, for our purposes, is that there is a $W$-associahedron for each finite Coxeter group $W$ and that the $W$-associahedron for $W=S_{n+1}$ is the usual $n$-dimensional associahedron.
Since we also constructed the usual associahedron from a lattice congruence of the weak order on the symmetric group, and since there is a weak order on every Coxeter group, it is natural to wonder whether we can generalize the lattice-theoretic construction.
That is, we would like a construction that, starting with an arbitrary Coxeter group $W$, produces a lattice congruence of the weak order on $W$ such that the Hasse diagram of the quotient is the $1$-skeleton of the generalized associahedron for $W$.

In this section, after filling in some background material on Coxeter groups, we describe the construction of the desired congruence, which we call a \emph{Cambrian congruence}.
We then discuss some properties of the Cambrian congruences and their quotients, the \emph{Cambrian lattices}.
Our coverage of Coxeter groups is necessarily no more than a sketch.
A more in-depth, but still gentle, introduction, including the definition of the generalized associahedron, is found in \cite{rsga}.  
Hohlweg's chapter \cite{Hohlweg} in this volume is another accessible description of Coxeter groups, culminating in a different construction of generalized associahedra.

To define a Coxeter group, choose a finite generating set $S=\set{s_1,\ldots,s_n}$ and for every $i<j$, choose an integer $m(i,j)\ge 2$, or $m(i,j)=\infty$.
Define $W$ to be the group with the following presentation:
\[W=\br{S\,\,|\,\,s_i^2=1,\,\,\forall \,i\,\,\mbox{ and }\,\,(s_is_j)^{m(i,j)}=1,\,\,\forall\, i<j}.\]
By convention, the ``relation'' $(s_is_j)^\infty=1$ is interpreted as the absence of a relation of the form $(s_is_j)^m=1$.
This abstract definition captures the essence of groups of transformation generated by reflections by taking a generating set of involutions and requiring that the product of two generators resembles the composition of two reflections (a rotation).

We mentioned above that the symmetric group is a Coxeter group.
Specifically, we take $S=\set{s_1,\cdots s_n}$, where  $s_i$ is the transposition $(i\,\,\,i\!+\!1)$.
We easily compute the order of each $s_is_j$ and conclude that we must set 
\[m(i,j)=\left\lbrace\begin{array}{ll}
3&\mbox{if }j=i+1,\mbox{ or}\\
2&\mbox{if }j>i+1.
\end{array}\right.\]
One can check that $S_{n+1}$ is indeed isomorphic to the abstract Coxeter group with this choice of $S$ and $m(i,j)$.

Besides $S_{n+1}$, we will follow another example of a Coxeter group:  The dihedral group of order $8$, which we call $B_2$.
This Coxeter group has $S=\set{s_1,s_2}$ and $m(1,2)=4$.
Thus $B_2$ is $\br{\set{s_1,s_2}\,\,|\,\,s_1^2=s_2^2=(s_1s_2)^4=1}$, and its elements are
\[1,\,s_1,\, s_2,\, s_1s_2,\,s_2s_1,\,s_1s_2s_1,\,s_2s_1s_2,\,s_1s_2s_1s_2=s_2s_1s_2s_1.\]

The \emph{Coxeter diagram} of $W$ is a graph with vertex set $\set{1,\ldots,n}$ and edges $i$ --- $j$ whenever $m(i,j)\ge 3$.
Each edge is labeled by $m(i,j)$, except that, by convention, we omit edge labels ``$3$.''
Note that an edge between $i$ and $j$ is absent if and only if $m(i,j)=2$, which happens if and only if $s_i$ and $s_j$ commute.
For example, the dihedral group of order $8$ has a diagram with two vertices connected by an edge labeled 4, and the diagram for $S_{n+1}$ is 
\[\scalebox{.8}{\includegraphics{diagramAn.ps}}\]

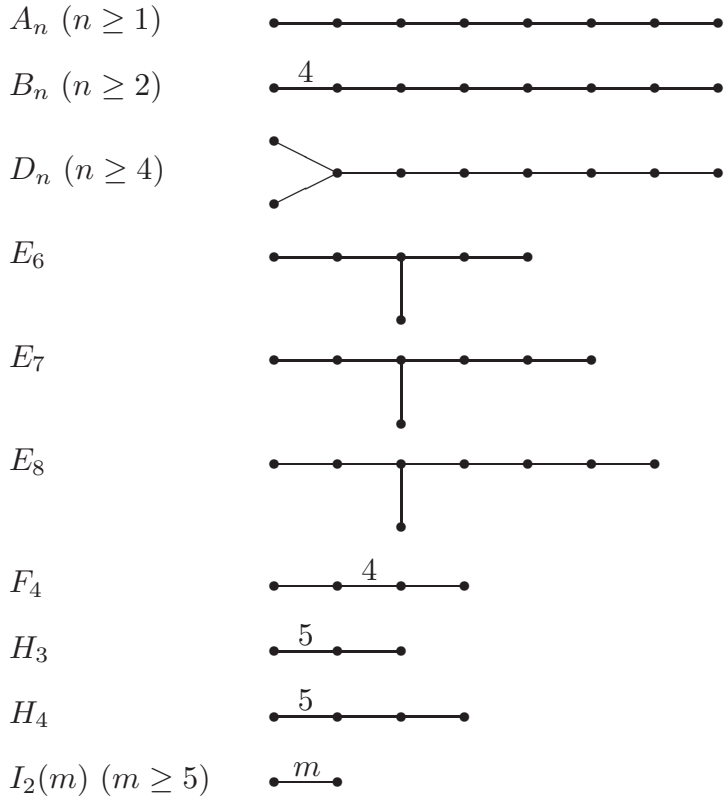
\begin{figure}
\begin{center}
\vspace{-.2in} 
\[ 
\begin{array}{lcl} 
A_n\  (n\geq 1) && 
\setlength{\unitlength}{1.00pt} 
\begin{picture}(140,17)(0,-2) 
\put(0,-0.15){\line(1,0){140}} 
\multiput(0,0)(20,0){8}{\circle*{3}} 
\end{picture}\\
B_n\ (n\geq 2)
&& 
\setlength{\unitlength}{1.00pt} 
\begin{picture}(140,17)(0,-2) 
\put(0,-0.15){\line(1,0){140}}
\multiput(0,0)(20,0){8}{\circle*{3}} 
\put(7.5,2){$4$}
\end{picture} \\[6 pt] 
D_n\ (n\geq 4)
&& 
\setlength{\unitlength}{1.00pt} 
\begin{picture}(140,17)(0,-2) 
\put(20,-0.15){\line(1,0){120}} 
\put(0,10){\line(2,-1){20}} 
\put(0,-10){\line(2,1){20}} 
\multiput(20,0)(20,0){7}{\circle*{3}} 
\put(0,10){\circle*{3}} 
\put(0,-10){\circle*{3}} 
\end{picture} 
\\[6 pt] 
E_6 
&& 
\setlength{\unitlength}{1.00pt} 
\begin{picture}(140,17)(0,-2) 
\put(0,-0.15){\line(1,0){80}} 
\put(40,0){\line(0,-1){20}} 
\put(40,-20){\circle*{3}} 
\multiput(0,0)(20,0){5}{\circle*{3}} 
\end{picture} 
\\[12 pt] 
E_7 
&& 
\setlength{\unitlength}{1.00pt} 
\begin{picture}(140,17)(0,-2) 
\put(0,-0.15){\line(1,0){100}} 
\put(40,0){\line(0,-1){20}} 
\put(40,-20){\circle*{3}} 
\multiput(0,0)(20,0){6}{\circle*{3}} 
\end{picture} 
\\[12pt] 
E_8 
&& 
\setlength{\unitlength}{1.00pt} 
\begin{picture}(140,17)(0,-2) 
\put(0,-0.15){\line(1,0){120}} 
\put(40,0){\line(0,-1){20}} 
\put(40,-20){\circle*{3}} 
\multiput(0,0)(20,0){7}{\circle*{3}} 
\end{picture} 
\\[18 pt] 
F_4 
&& 
\setlength{\unitlength}{1.00pt} 
\begin{picture}(140,17)(0,-2) 
\put(0,-0.15){\line(1,0){60}} 
\multiput(0,0)(20,0){4}{\circle*{3}} 
\put(27.5,2){$4$}
\end{picture} 
\\
H_3
&&
\setlength{\unitlength}{1.00pt} 
\begin{picture}(140,17)(0,-2) 
\put(0,-0.15){\line(1,0){40}} 
\multiput(0,0)(20,0){3}{\circle*{3}} 
\put(7.5,2){$5$}
\end{picture} 
\\
H_4
&&
\setlength{\unitlength}{1.00pt} 
\begin{picture}(140,17)(0,-2) 
\put(0,-0.15){\line(1,0){60}} 
\multiput(0,0)(20,0){4}{\circle*{3}} 
\put(7.5,2){$5$}
\end{picture} 
\\
I_2(m)\  (m\geq 5) 
&&
\setlength{\unitlength}{1.00pt} 
\begin{picture}(140,17)(0,-2) 
\put(0,-0.15){\line(1,0){20}} 
\multiput(0,0)(20,0){2}{\circle*{3}} 
\put(6,2){$m$}
\end{picture} 
\end{array} 
\] 
\end{center}
\vspace{-.1in} 
\caption{Coxeter diagrams of finite irreducible Coxeter systems}
\label{Cox diagrams} 
\end{figure}

Coxeter diagrams provide a convenient way of encoding the defining data of a Coxeter group.
This encoding makes it easy to describe one of the fundamental results on Coxeter groups:  The classification of finite Coxeter groups.
If the Coxeter diagram of $W$ is not connected as a graph, then it is easy to see that $W$ is a direct product of the Coxeter groups encoded by the connected components of the diagram.
If the Coxeter diagram of $W$ is connected, then $W$ is called irreducible.
Figure~\ref{Cox diagrams} lists the Coxeter diagrams of finite irreducible Coxeter groups and shows their standard names.
In particular, the name $B_2$ for our earlier example is part of this naming convention, and the Coxeter group $S_{n+1}$ appears as $A_n$ in the classification.

\begin{figure}
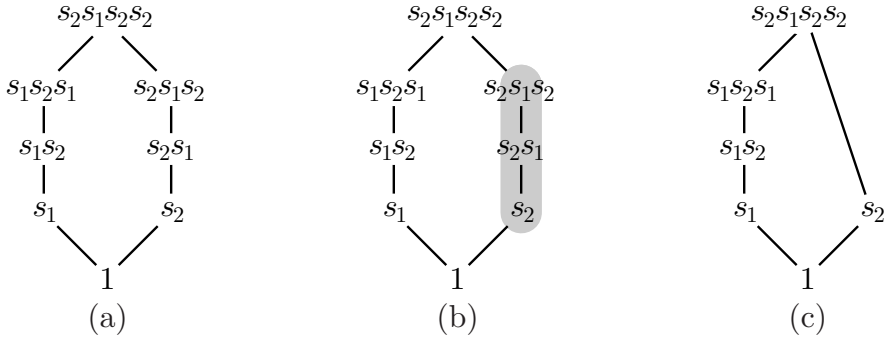

\begin{tabular}{ccc}
\scalebox{1}{\includegraphics{weakB2.ps}}&
\scalebox{1}{\includegraphics{weakB2cong.ps}}&
\scalebox{1}{\includegraphics{B2camb.ps}}\\
(a)&(b)&(c)
\end{tabular}
\caption{a:  The weak order on $B_2$.  b:  A Cambrian congruence on $B_2$.  c:  The corresponding Cambrian lattice.}
\label{weakB2 fig}
\end{figure}

Since $W$ is generated by $S$, each element $w$ of $W$ can be written (in many ways) as a word in the ``alphabet'' $S$.
A word of minimal length, among words for $w$, is called a \emph{reduced word} for $w$.
The \emph{length} $\ell(w)$ of $w$ is the length of a reduced word for $w$.
The \emph{weak order} on a Coxeter group $W$ sets $u\le w$ if and only if a reduced word for $u$ occurs as a prefix of some reduced word for $w$. 
The cover relations in the weak order are $w\covered ws$ for $w\in W$ and $s\in S$ with $\ell(w)<\ell(ws)$.
It is easy to check that this definition of the weak order reduces, in the case where $W=S_{n+1}$, to the earlier definition of the weak order on permutations.
(See Figure~\ref{tam-weakS4 fig}.b.)
The weak order on $B_2$ is pictured in Figure~\ref{weakB2 fig}.
To see that this picture is correct, keep in mind that the element $s_1s_2s_1s_2$ at the top of the picture is also equal to $s_2s_1s_2s_1$.
The weak order is a meet semilattice (meaning that meets exist, but not necessarily joins) in general, and a lattice when $W$ is finite.
(What we have defined is sometimes called the \emph{right weak order}.  
There is also an isomorphic, but not identical, \emph{left weak order}.
See, for example, Forcey's chapter \cite{Forcey} in this volume.)

\begin{figure}
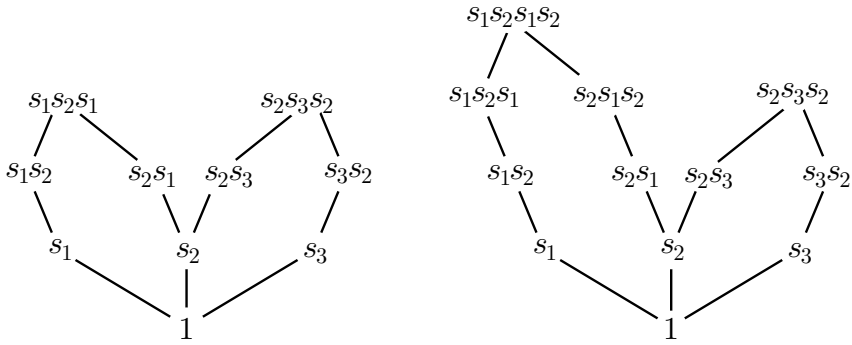

\centering
\scalebox{1}{\includegraphics{weakS4diagram.ps}}
\qquad
\scalebox{1}{\includegraphics{weakB3diagram.ps}}
\caption{The Coxeter diagram inside the weak order.  On the left, the Coxeter group is $A_3=S_4$ and on the right is $B_3$.}
\label{weakdiagram fig}
\end{figure}

Having defined Coxeter diagrams and the weak order, we are prepared to state the general definition of Cambrian lattices.
Informally:  we orient the Coxeter diagram, find a copy of the Coxeter diagram within the weak order, and contract edges according to the orientation.
More formally, notice that for each pair $s_i,s_j$ of generators in $S$, the datum $m(i,j)$ is easily read off from the weak order by finding the join (least upper bound) of $s_i$ and $s_j$.
The join $s_i\join s_j$ is of the form $s_is_js_is_j\cdots$ and has length $m(i,j)$.
The interval in the weak order below $s_i\join s_j$ has the form of a polygon with $2m(i,j)$ sides.
The union of all of these intervals, as $i$ and $j$ vary, is essentially the Coxeter diagram.
To make the resemblance more perfect, we leave out those intervals where $m(i,j)=2$, or in other words, those intervals that form $4$-gons.
For example, Figure~\ref{weakdiagram fig} shows this union of intervals in $A_3=S_4$ and in $B_3$.

Now recall the characterization, from Section~\ref{cong sec}, of the Tamari congruence on $S_{n+1}$ as the smallest congruence contracting a certain set of edges.
These edges are $s_2\covered s_2s_1$, $s_3\covered s_3s_2$, \ldots, $s_n\covered s_ns_{n-1}$.
More generally, as we let the polygon $Q$ vary in the permutations-to-triangulations map, there is a choice of $s_1\covered s_1s_2$ or $s_2\covered s_2s_1$, $s_2\covered s_2s_3$ or $s_3\covered s_3s_2$, and so forth, such that the fibers of $\eta$ constitute the smallest congruence contracting the chosen edges.

Accordingly, we define the \emph{Cambrian congruences} on a general finite Coxeter group $W$ as follows:
For each edge $i$------$j$ in the Coxeter diagram of $W$, declare one of $i$ or $j$ to be ``before'' the other.
We may as well name $i$ and $j$ so that $i$ is before $j$.
There is a chain $s_j\covered s_js_i\covered\cdots\covered (s_js_is_js_i\cdots)$ in the weak order with the top element having $m(i,j)-1$ letters.
We contract all covers in this chain.
After making such a choice for each edge of the diagram, the Cambrian congruence is the finest congruence of the weak order on $W$ contacting all of the chosen edges in the weak order.
In $S_{n+1}$, edges in the diagram connect $(i-1)$ to $i$ for $i=2,\ldots,n$.
We declare $i-1$ to be before $i$ if $i$ is lower-barred and declare $i$ to be before $i-1$ if $i$ is upper-barred.
Figure~\ref{weakB2 fig}.b shows one of the two Cambrian congruences on $B_2$.
The contracted edges are indicated by shading.

\begin{figure}
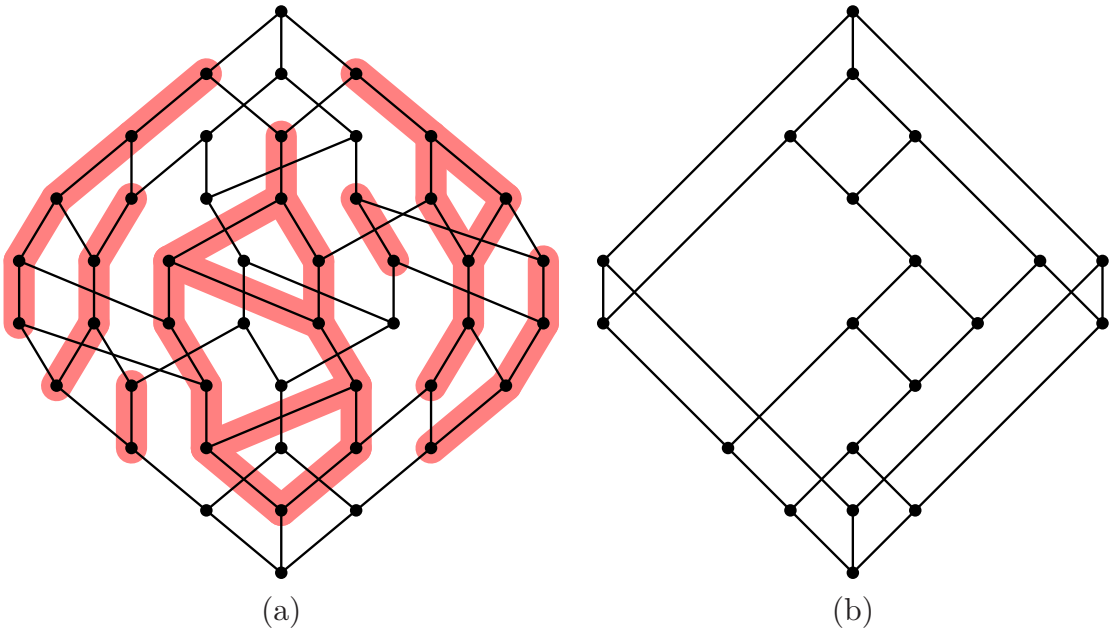

\begin{tabular}{cc}
\scalebox{.98}{\includegraphics{weakB3cong.ps}}&
\scalebox{.98}{\includegraphics{B3camb.ps}}\\
(a)&(b)
\end{tabular}
\caption{a:   A Cambrian congruence on $B_3$.  b:  The corresponding Cambrian lattice.}
\label{B3camb fig}
\end{figure}

To understand larger examples beyond the $S_{n+1}$ case, we need to know how edge-forcing works for the weak order on a general finite Coxeter group.
In general, we must consider all \textit{polygonal intervals} (intervals whose Hasse diagrams are cycles) in the weak order.
The weak order on permutations has the property that every polygonal interval is a quadrilateral or a hexagon. 
In a general finite Coxeter group, larger polygonal intervals may be present.
Edge forcing in general is completely analogous to edge-forcing for the weak order on permutations.
There are local forcing requirements for every polygonal interval, with the same description: 
Whenever a bottom or top edge of the interval is contracted, the opposite bottom or top edge is also contracted, as well as all side edges (if side edges are present).
The only difference is that, in general, there might be more than two side edges.
We define \emph{closure under local forcing} in the same way, with reference to all polygonal intervals, and again, \cite[Theorem~25]{hyperplane} implies that a collection $\E$ of edges in the weak order on $W$ is the set of edges contracted by some congruence if and only if $\E$ is closed under local forcing.
Figure~\ref{B3camb fig}.a shows the weak order on the Coxeter group $B_3$, with certain edges shaded.
Using the local forcing criterion and comparing to Figure~\ref{weakdiagram fig}, the reader can readily verify that the congruence indicated by the shaded edges is the Cambrian congruence given by orienting $1$ before $2$ and $3$ before~$2$.

A \emph{Cambrian lattice} is the quotient of the weak order on a finite Coxeter group $W$ modulo some Cambrian congruence on $W$.
A general fact about quotients of finite lattices says that this quotient is isomorphic to the subposet induced by the bottom elements of congruence classes.
Above, we described how to define a Cambrian congruence on $S_{n+1}$ starting from a choice of the polygon $Q$.
The Cambrian lattice obtained from this orientation is isomorphic to the Tamari-like lattice on triangulations of $Q$.
Thus two Cambrian lattices associated to $S_4$ are shown in Figures~\ref{tam-weakS4 fig}.a and~\ref{permtri2-nontam fig}.b.
As further examples, Figure~\ref{weakB2 fig}.c depicts the Cambrian lattice arising from the Cambrian congruence shown in Figure~\ref{weakB2 fig}.b, while Figure~\ref{B3camb fig}.b shows the Cambrian lattice arising from the Cambrian congruence shown in Figure~\ref{B3camb fig}.a.

It should seem unlikely, \emph{a priori}, that this generalization of the Tamari lattice along lattice-theoretic lines should yield anything useful.  
But amazingly, the Cambrian lattices have the same relationship to the generalized associahedra that the Tamari lattice has to the usual associahedron.
The following theorem was conjectured in \cite{cambrian} and proved in \cite{sortable,sort_camb,camb_fan}.
\begin{theorem}\label{camb Hasse}
The Hasse diagram of any Cambrian lattice associated to a finite Coxeter group $W$ is isomorphic to the $1$-skeleton of the generalized associahedron for $W$.
\end{theorem}


\section{Cambrian fans}\label{fan sec}
Theorem~\ref{camb Hasse} makes the connection, in general, between Cambrian lattices and generalized associahedra.
In this section, we expand on that connection by describing how a Cambrian congruence determines a polyhedral structure closely related to the generalized associahedron:  a Cambrian fan.

Consider a finite Coxeter group $W$, represented concretely as a group generated by the reflections $S$ in a real vector space of dimension $|S|$.
Typically, there are additional elements of $W$, besides the generators $S$, that act as reflections.
Let $T$ be the set of elements of $W$ that act as reflections.
For each $t\in T$, there is a corresponding reflecting hyperplane. 
Let $\A$ be the collection of all of these reflecting  hyperplanes. 
This is called the \emph{Coxeter arrangement} associated to $W$.

\begin{figure}
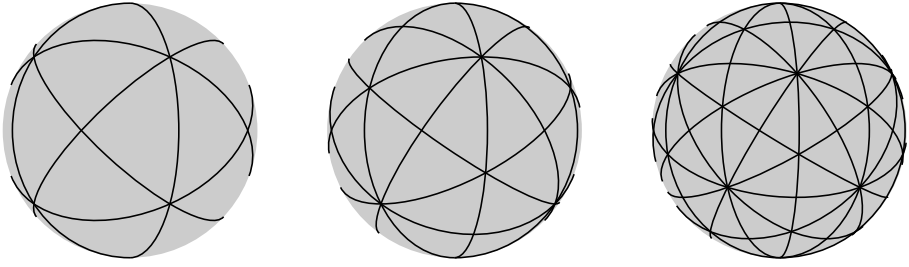

\centering
\scalebox{1}{\includegraphics{a3.ps}}\qquad\scalebox{1}{\includegraphics{b3.ps}}\qquad\scalebox{1}{\includegraphics{h3.ps}}
\caption{The Coxeter arrangements for $A_3$, $B_3$, and $H_3$}
\label{Coxarr fig}
\end{figure}

Referring to Figure~\ref{Cox diagrams}, we see that there are three irreducible Coxeter groups with $|S|=3$.
Figure~\ref{Coxarr fig} pictures the associated Coxeter arrangements.
Each arrangement is a collection of planes, through the origin, in $\reals^3$.
The intersection of these planes with the unit sphere about the origin is an arrangement of great circles.
The figure shows these great circles on the sphere.
The sphere is opaque, so that we only see each great circle as it intersects the near side of the sphere.

The complement $\reals^n\setminus\left(\bigcup_{H\in\A}H\right)$ of the arrangement $\A$ is a collection of unbounded, $n$-dimensional, disjoint, open sets.
The closures of these sets are called \emph{regions}.
The regions are the maximal cones of a \emph{fan}, meaning that any two regions intersect in a \emph{face} of each.
We call this fan the \emph{Coxeter fan} $\F$ associated to $W$.
In Figure~\ref{Coxarr fig}, the regions appear as spherical triangles, each representing an unbounded cone with triangular cross section.
In general, the regions are unbounded cones whose cross sections are $(n-1)$-dimensional simplices, and therefore $\F$ is a \emph{simplicial fan}.

The hyperplane arrangement $\A$ and the fan $\F$ it defines are closely related to the combinatorics of the Coxeter group $W$.
We now highlight two well-known aspects of this close relation.
First and most importantly, the regions are in bijection with the elements of $W$.
There is a region $D$ whose facets are contained in the reflecting hyperplanes for the reflections $S$.
(For some Coxeter groups $W$, it is possible to choose a representation of $W$ such that no such $D$ exists, but for any $W$, there is a standard way of producing a representation of $W$ as a group generated by reflections, and in this standard representation, $D$ exists.)
Then each element $w$ in $W$ is mapped to the region $wD$ obtained by acting on $D$ by the transformation $w$.

\begin{figure}
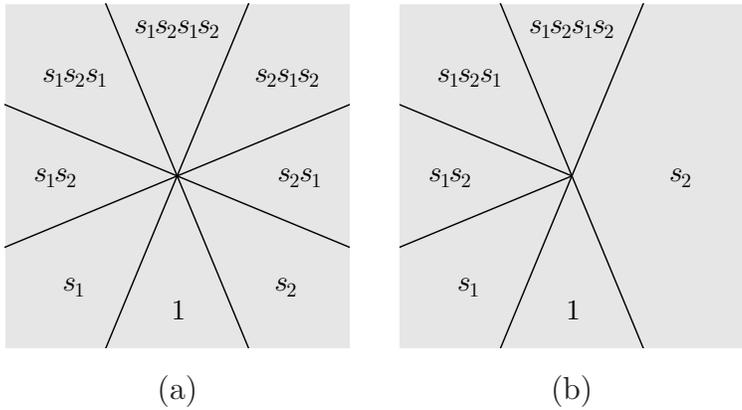

\centerline{
\begin{tabular}{cc}
\scalebox{.9}{
\includegraphics{B2F.ps}
}
&
\scalebox{.9}{
\includegraphics{B2Fc.ps}
}
\\[4 pt]
(a)&(b)
\end{tabular}
}
\caption{a:  The Coxeter fan for $B_2$.  b:  A Cambrian fan associated to $B_2$}
\label{B2Fc fig}
\end{figure}

\begin{figure}
\scalebox{.55}{\includegraphics{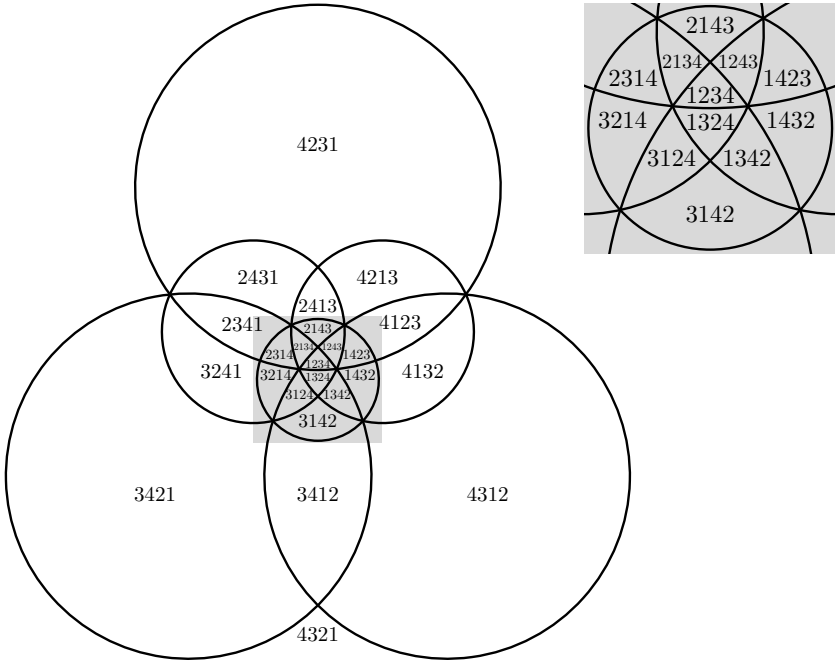}}
\caption{Regions in the Coxeter arrangement for $S_4$, labeled by permutations.}
\label{A3labels fig}
\end{figure}

As an example, consider the case where $W$ is $B_2$.
The Coxeter fan for $B_2$ is shown in Figure~\ref{B2Fc fig}.a.
As another example, we give a different depiction of the Coxeter arrangement for $A_3$, which appeared as the left picture in Figure~\ref{Coxarr fig}.
As before, we first pass from a collection of six planes through the origin in $\reals^3$ to a collection of six great circles on the unit sphere.
But in Figure~\ref{A3labels fig}, instead of showing only the near side of an opaque sphere, we show a stereographic projection of the sphere to the plane. 
This allows us to label the 24 regions (including the region shown as an unbounded region outside of all of the circles shown) by the 24 permutations in $S_4$.
Since some labels are small, in the top-right corner of the figure we show a magnified version of the center of the picture.

The second aspect we wish to highlight is that cover relations in the weak order correspond to pairs $Q,R$ of adjacent regions.
The pair $Q,R$ is separated by a unique hyperplane $H\in\A$.
Taking, without loss of generality, $Q$ to be on the same side of $H$ as $D$, the cover relation is $Q\covered R$.
This is readily verified in our running examples by comparing Figure~\ref{B2Fc fig}.a to Figure~\ref{weakB2 fig} and comparing Figure~\ref{A3labels fig} to Figure~\ref{tam-weakS4 fig}.b.

\begin{figure}
\scalebox{.55}{\includegraphics{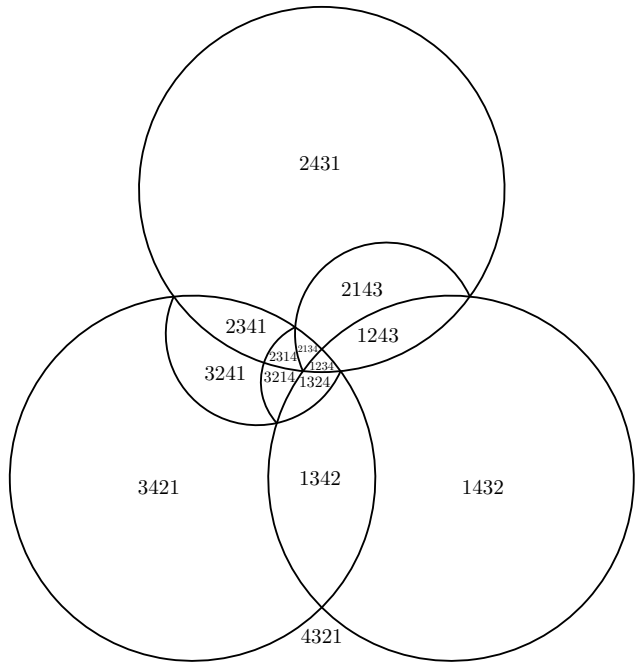}}
\caption{The $s_1s_2s_3$-Cambrian fan associated to $S_4$}
\label{A3camb fig}
\end{figure}

Lattice congruences on the weak order respect the geometry of $\F$.
Recall that every congruence class is an interval.
Using the geometric characterization of the weak order, it is not hard to show that, given any interval in the weak order, the union of the corresponding regions in $\F$ is a convex cone.
Thus the lattice congruence defines a collection of convex full-dimensional cones, and \cite[Theorem~1.1]{con_app} says that these are the maximal cones of a fan.
When the lattice congruence is a Cambrian congruence, the resulting fan is called a \emph{Cambrian fan}.
Figure~\ref{B2Fc fig}.b illustrates the Cambrian fan associated to the Cambrian congruence on $B_2$ shown in Figure~\ref{weakB2 fig}.b.
In the figure, each maximal cone is labeled by the bottom element of the corresponding Cambrian congruence class.
A Cambrian fan associated to $S_4$ is shown in Figure~\ref{A3camb fig}.
This is the fan associated to the Cambrian congruence (or Tamari congruence) pictured in Figure~\ref{permtri-tamcong fig}.b.
The fan is drawn in stereographic projection as explained in connection with Figure~\ref{A3labels fig}, and each maximal cone of the fan is labeled with the bottom element of the corresponding congruence class.

\begin{figure}
\scalebox{.7}{\includegraphics{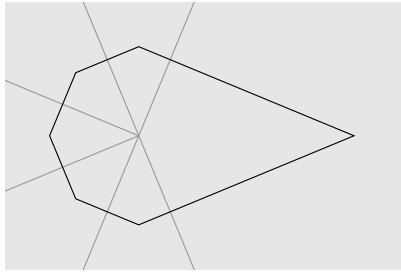}}
\caption{A $B_2$-associahedron and its normal fan (a Cambrian fan)} 
\label{B2assocnorm fig}
\end{figure}

General considerations about fans constructed from lattice congruences of the weak order imply the following connection between Cambrian fans and Cambrian lattices:
Two maximal cones in the Cambrian fan are adjacent (share a codimension-$1$ face) if and only if the corresponding elements of the Cambrian lattice are related by a cover.
In particular, by Theorem~\ref{camb Hasse}, the adjacency graph on maximal cones in the Cambrian fan is isomorphic to the $1$-skeleton of the generalized associahedron.
In fact, more is true, as conjectured in \cite{cambrian} and proved in \cite{sortable,sort_camb,camb_fan}:
\begin{theorem}\label{camb_fan thm}
Each Cambrian fan associated to $W$ is combinatorially isomorphic to the normal fan of the $W$-associahedron.
\end{theorem}
Rather than actually defining the normal fan of a polytope, let us informally describe the relationship between a polytope $P$ and its normal fan $\F$.
First, there is a one-to-one correspondence between facets (maximal proper faces) of $P$ and rays in the fan $\F$ such that each ray $\rho$ consists of outward-facing normal vectors to the corresponding facet $F$.
Then, there is a one-to-one correspondence between codimension-$2$ faces of $P$ and two-dimensional cones in $\F$ such that a ray is contained in a two-dimensional cone if and only if the corresponding facet contains the corresponding codimension-$2$ face.
Continuing, similar correspondences and reversals of containment are required in each dimension.
We illustrate in Figure~\ref{B2assocnorm fig} by picturing a realization of the $B_2$-associahedron whose normal fan is the Cambrian fan pictured in Figure~\ref{B2Fc fig}.b.

The combinatorial isomorphism with the normal fan of the $W$-associahedron is not the ultimate result on Cambrian fans. 
Indeed, there is a natural polytopal realization of the $W$-associahedron whose normal fan actually coincides with the Cambrian fan.
Hohlweg's chapter \cite{Hohlweg} in this volume describes this realization.

\section{Sortable elements}\label{sort sec}
Although the definition of Cambrian lattices is a lattice-theoretic one, most of what has been proved about them follows from combinatorial models.
In the paper~\cite{cambrian}, where Cambrian lattices were first defined, the Cambrian lattices associated to the Coxeter groups $S_{n+1}=A_n$ and $B_n$ were described in terms of the combinatorics of triangulations, with the help of some results from the theory of fiber polytopes \cite{Fiber,Iterated,Equivariant}.
A general model for Cambrian lattices was later provided by the \emph{sortable elements}.
These were introduced in~\cite{sortable}, and the connection to Cambrian lattices was made in~\cite{sort_camb}.

In defining a Cambrian lattice, our initial datum is an orientation of the Coxeter diagram of $W$.
To define sortable elements, it is useful to write this initial datum in a different form.
A \emph{Coxeter element} is an element $c$ of $W$ that has a reduced word of the form $s_1\cdots s_n$, where $S=\set{s_1,\ldots,s_n}$ and $|S|=n$.
That is, a Coxeter element is the product of the generators in $S$, each occurring exactly once.
Choosing a Coxeter element $c$ in a finite Coxeter group $W$ is equivalent to choosing, for each non-commuting pair of generators $s_i$ and $s_j$, which of the two is before the other in a word for $c$.
Since $s_i$ and $s_j$ are non-commuting if and only $i$---$j$ is an edge in the Coxeter diagram, Coxeter elements are equivalent to orientations of the diagram.
For this reason, we use the term \emph{$c$-Cambrian congruence} to describe the Cambrian congruence arising from the diagram orientation corresponding to $c$.
As an example, if we orient the diagram $s_1$------$s_2$------$s_3$ of $S_4$ with $s_1$ before $s_2$ and $s_3$ before $s_2$, we define the Coxeter element $c=s_1s_3s_2=s_3s_1s_2$.

In the case of the symmetric group, the Coxeter element $c$ can be read off directly from the polygon $Q$.
Specifically, $c$ is the $(n+1)$-cycle in $S_{n+1}$ obtained by reading the vertices of $Q$ (excluding $0$ and $n+2$) in counter-clockwise order.

To define sortable elements, it is useful to think about how one might write down a reduced word for an element $w$ of $W$.
We assume that we can easily determine the length of a given element of $W$.
(Although this assumption seems to ignore the fact that the length of $w$ is defined in terms of reduced words, it is nevertheless a useful point of view.)
If $w$ is the identity, then the only reduced word for $w$ is the empty word.
Otherwise, we write down a reduced word $a_1\cdots a_k$ for $w$ from left to right as follows.
First, try each element $s_i\in S$ in some order until we find one such that $\ell(s_iw)<\ell(w)$.
We will eventually find such an $s_i$ because $w$ has some nonempty reduced word, and we can take $s_i$ to be the first letter of that word.
Set $a_1$ equal to $s_i$, and write $w'=a_1w$.
If $w'$ is the identity, then $k=1$ and $a_1$ is the desired reduced word.
Otherwise, find $s_{i'}$ such that $\ell(s_{i'}w')<\ell(w')$, set $a_2$ equal to $s_{i'}$, and define $w''=a_2w'$.
If $w''$ is the identity, then $k=2$ and $a_1a_2$ is the desired word.
Otherwise, continuing in this manner, we eventually find a reduced word $a_1\cdots a_k$.
The output depends, of course, on the order in which we try the elements of $S$ in each step. 

The usefulness of this method of writing down a reduced word is that it lends itself to a global choice of a canonical reduced word for each element.
For example, one might, at each step, try the elements of $S$ in the order $s_1,\ldots,s_n$.
This would have the effect of representing each element by its lexicographically first reduced word, in the sense of lexicographic order on subscripts.

The definition of sortable elements depends on a more subtle choice of a canonical word.
Fix a total order $(s_1,\ldots,s_n)$ on $S$.
Follow the method described above, trying the letters in $S$ cyclically in the order $(s_1,\ldots,s_n)$.
This method differs from the lexicographic method described above, in that we do not start again at $s_1$ at every step.
Instead, if we choose $s_i$ at some step, in the next step we try the letters in the order $s_{i+1},\ldots,s_n,s_1,\ldots,s_i$.
The resulting reduced word is called the \emph{$(s_1,\ldots,s_n)$-sorting word} for $w$.
In effect, we take repeated ``passes'' through the sequence $(s_1,\ldots,s_n)$, adding whatever letters we can to the word.
For convenience, between passes, we insert the symbol ``\,$|$\,'' as a ``divider.'' 

An $(s_1,\ldots,s_n)$-sorting word can be described by a sequence of subsets of $S$:
The $i\th$ subset in the sequence is the set of letters that are added to the $(s_1,\ldots,s_n)$-sorting word in the $i\th$ pass.
Equivalently, these are the sets of letters occurring between dividers.
Let $c$ be the Coxeter element $s_1\cdots s_n$.
We say that $w$ is \emph{$c$-sortable} if this sequence of subsets is weakly decreasing under containment.
(We say ``$c$-sortable'' instead of ``$(s_1,\ldots,s_n)$-sortable'' because it is easy to show that the notion depends only on $c$, not on the reduced word $s_1\cdots s_n$ for $c$.)

\begin{figure}
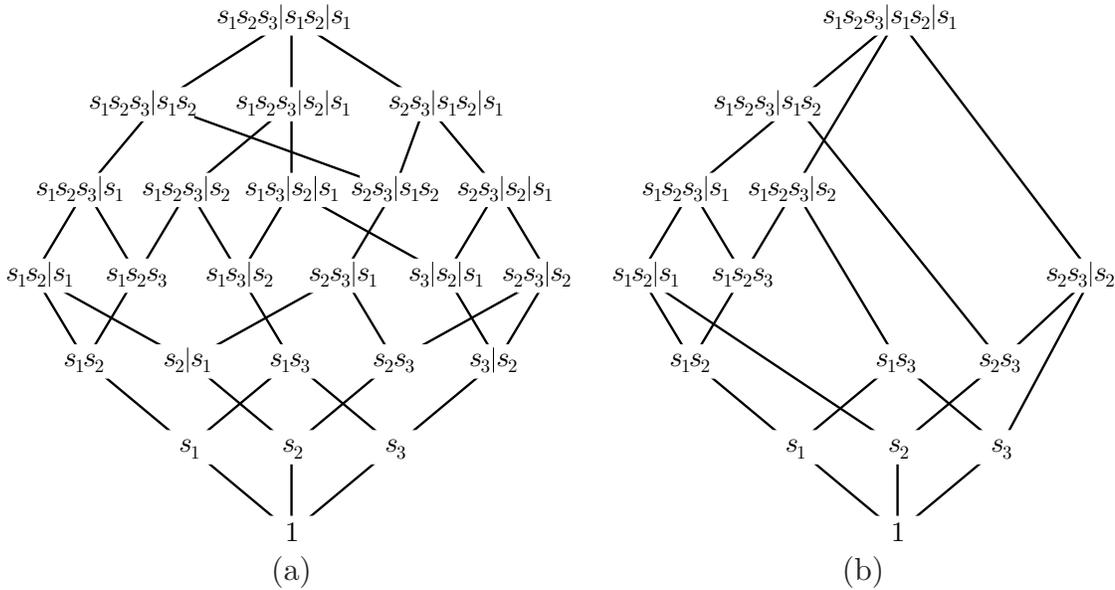

\begin{tabular}{cc}
\scalebox{.98}{\includegraphics{weakS4sorting_abc.ps}}&
\scalebox{.98}{\includegraphics{weakS4bottoms_abc.ps}}\\
(a)&(b)
\end{tabular}
\caption{a:  Sorting words for elements of $S_4$.  b: The subposet of the weak order induced by $s_1s_2s_3$-sortable elements.}
\label{weakS4sorting fig}
\end{figure}

The $(s_1,s_2)$-sorting words for elements of $B_2$ are the empty word, $s_1$, $s_1s_2$, $s_1s_2|s_1$, $s_1s_2|s_1s_2$, $s_2$, $s_2|s_1$, and $s_2|s_1s_2$.
Six of the elements of $B_2$ are $c$-sortable for $c=s_1s_2$.
The exceptions are $s_2|s_1$ (because $\set{s_2}\not\supseteq\set{s_1}$) and $s_2|s_1s_2$ (because $\set{s_2}\not\supseteq\set{s_1,s_2}$).
The $(s_1,s_2,s_3)$-sorting words for elements of $S_4$ are shown in Figure~\ref{weakS4sorting fig}.a, in the same positions as the permutations shown in Figure~\ref{tam-weakS4 fig}.b.
Figure~\ref{weakS4sorting fig}.b shows the subposet of the weak order induced by the $s_1s_2s_3$-sortable elements.

Comparing Figure~\ref{weakS4sorting fig}.a with Figure~\ref{permtri-tamcong fig}.b, we see an example of the following theorem, which is the concatenation of \cite[Theorem~1.1]{sort_camb} and \cite[Theorem~1.4]{sort_camb}.
\begin{theorem}\label{sort camb thm}
An element $w$ of $W$ is the bottom element in its $c$-Cambrian congruence class if and only if it is $c$-sortable.
\end{theorem}
As mentioned earlier, the quotient of a finite lattice modulo a congruence $\Theta$ is isomorphic to the subposet induced by the bottom elements of $\Theta$-classes.
Thus, we have the following corollary, which is exemplified by Figures~\ref{tam-weakS4 fig}.a and~\ref{weakS4sorting fig}.b.
\begin{cor}\label{sort camb induced}
The $c$-Cambrian lattice is isomorphic to the subposet of the weak order on $W$ induced by the $c$-sortable elements.
\end{cor}
Typically, the bottom elements of a lattice congruence on a finite lattice $L$ are not a sublattice of $L$.
But \cite[Theorem~1.2]{sort_camb} shows that the Cambrian congruence is an exception:
\begin{theorem}\label{camb sub thm}
The $c$-sortable elements in $W$ constitute a sublattice of the weak order on $W$.
\end{theorem}

The realization of the Cambrian lattice in terms of sortable elements provides considerable combinatorial traction.
Most importantly, the sortable elements allow for \emph{uniform} proofs.
A non-uniform (or \emph{case-by-case} or \emph{type-by-type}) proof of a result on finite Coxeter groups is a proof obtained by checking each case of the classification of finite Coxeter groups (Figure~\ref{Cox diagrams}).
Typically, the result is argued separately for the infinite families (in order of difficulty) $A_n$, $B_n$ and $D_n$, and then separately for the exceptional cases $E_6$, $E_7$, $E_8$, $F_4$, $H_3$, $H_4$ and $I_2(m)$.
A uniform proof is a single argument that is valid for an arbitrary finite Coxeter group.

We now describe in detail three specific types of combinatorial traction provided by sortable elements.
The first is a natural recursive structure that lends itself well to arguments by induction on two parameters:  the length of elements of $W$ and the rank of $W$.
(The rank of a Coxeter group is the size of its defining generating set $S$.)
The recursive structure is encapsulated in two easy lemmas \cite[Lemmas~2.4,~2.5]{sortable}, which are given below as Lemmas~\ref{sc} and~\ref{scs}. 

Suppose $s\in S$ and let $c$ be a Coxeter element of $W$.
Then $s$ is \emph{initial} in $c$ if there is a reduced word for $c$ having $s$ as its first letter.
Equivalently, in the diagram orientation corresponding to $c$, no other generator is declared to be ``before'' $s$.
The notation $W_{\br{s}}$ stands for the subgroup of $W$ generated by the set $S\setminus\set{s}$.
This is an example of a \emph{standard parabolic subgroup} of $W$, and it is a Coxeter group in its own right, with generating set $S\setminus\set{s}$.

\begin{lemma}
\label{sc}
Let~$s$ be initial in~$c$ and suppose $w\not\ge s$.
Then~$w$ is $c$-sortable if and only if it is an $sc$-sortable element of~$W_{\br{s}}$.
\end{lemma}

To see why this lemma is true, write $c=s_1\cdots s_n$ with $s_1=s$.
The assumption that $w\not\ge s$ means that there is no reduced word for $w$ having the word $s$ as a prefix, or in other words, there is no reduced word for $w$ having $s$ as its first letter.
In particular, the $(s_1,\ldots,s_n)$-sorting word for $w$ does not start with $s_1$.
If $w$ is $c$-sortable, then $s_1$ does not appear in its $(s_1,\ldots,s_n)$-sorting word, so $w$ is in $W_{\br{s}}$, and its $(s_2,\ldots,s_n)$-sorting word is the same sequence of letters as its $(s_1,\ldots,s_n)$-sorting word.
Now $sc=s_2\cdots s_n$, and the fact that $w$ is $c$-sortable immediately implies that $w$ is $sc$-sortable, as an element of $W_{\br{s}}$.
For the converse, we need a well-known fact about standard parabolic subgroups: 
If $w$ is contained in $W_{\br{s}}$, then no reduced word for $w$ (as an element of $W$) contains the letter $s$.
Thus if $w$ is in $W_{\br{s}}$ and is $sc$-sortable as an element of $W_{\br{s}}$, then its $(s_1,\ldots,s_n)$-sorting word must coincide with its $(s_2,\ldots,s_n)$-sorting word.
Now the $c$-sortability of $w$ follows from the $sc$-sortability of $w$.

\begin{lemma}
\label{scs}
Let~$s$ be initial in~$c$ and suppose $w\ge s$.
Then~$w$ is $c$-sortable if and only if $sw$ is $scs$-sortable.
\end{lemma}
To explain this lemma, again write $c=s_1\cdots s_n$ with $s_1=s$ and note that $scs$ is the Coxeter element $s_2\cdots s_ns_1$.
To make the $(s_1,\ldots,s_n)$-sorting word for $w$, we first need to try the letter $s_1$.
Under the hypotheses of the lemma, (with $s=s_1$), we take $s_1$ as the first letter of the $c$-sorting word, and then continue trying the letters in the cyclic order given by $(s_2,\ldots,s_n,s_1)$.
The sequence of letters obtained after $s_1$ is exactly the $(s_2,\ldots,s_n,s_1)$-sorting word for $sw$ and the criterion for $w$ to be $c$-sortable is exactly the criterion for $sw$ to be $scs$-sortable.

The second type of combinatorial traction provided by sortable elements is the natural \emph{search tree} structure on sortable elements.
Each $c$-sortable element has a naturally-defined predecessor among $c$-sortable elements, obtained by deleting the last letter of the $(s_1,\ldots,s_n)$-sorting word.
Starting with the Hasse diagram of the Tamari lattice, delete all edges except those connecting an element to its predecessor.
The result is a spanning tree of the Hasse diagram, rooted at the identity element.
This tree may depend on $(s_1,\ldots,s_n)$, not only on $c$.
Furthermore, given any $c$-sortable element $w$, there is a simple algorithm for determining all $c$-sortable elements whose predecessor is $w$.
Under these circumstances, it is a simple matter to efficiently traverse the set of $c$-sortable elements.
For a general description of efficient traversal, see \cite[Section~4C]{StembridgeMSJ}.
Figure~\ref{weakS4search fig} shows the search tree structure on $s_1s_2s_3$-sortable elements in $S_4$.
\begin{figure}
\scalebox{1}{\includegraphics{weakS4search_abc.ps}}
\caption{The search tree structure on $s_1s_2s_3$-sortable elements in $S_4$}
\label{weakS4search fig}
\end{figure}

The third type of combinatorial traction provided by sortable elements is the extra combinatorial information available in the sorting word for a $c$-sortable element.
We explain by giving an important example of information that can be read off from the sorting word.  
(Another example appears in Hohlweg's chapter \cite[Section~3.3]{Hohlweg}.)
To explain the example, we need to define the \emph{simple roots} associated to our chosen reflection representation of $W$.
For each generator $s_i\in S$, there is a corresponding vector $\alpha_i$ called a simple root. 
Recall that $D$ is the region identified with the identity element, and that its facets are contained in the reflecting hyperplanes for the generators $S$.
The full definition of a \emph{root system} specifies the length of each $\alpha_i$ up to a global scaling, but for our purposes, all we need to say is that $\alpha_i$ is a nonzero normal vector to the reflecting hyperplane for $s_i$, pointing in the direction of the interior of $D$.

Let $c=s_1\cdots s_n$ and let $v$ be a $c$-sortable element of $W$ with $(s_1,\ldots,s_n)$-sorting word $a_1\cdots a_k$.
Fix some $s_i\in S$.
In the process of building the $(s_1,\ldots,s_n)$-sorting word, there was some first time that $s_i$ was tested, but \textbf{not} included in the $(s_1,\ldots,s_n)$-sorting word.
Let $a_1\cdots a_j$ be the part of the $(s_1,\ldots,s_n)$-sorting word that was already determined before $s_i$ was first tested and not included.
We say $a_1\cdots a_k$ \emph{skips} $s_i$ after position $j$.
Define $C_c^{s_i}(v)$ to be the vector $a_1\cdots a_j \cdot \alpha_i$.
That is, act on the simple root $\alpha_i$ by the reflection $a_j$, then the reflection $a_{j-1}$, on so forth until acting by $a_1$.
As indicated by the notation, the vector $C_c^{s_i}$ depends only on $c$, not on $(s_1,\ldots,s_n)$.
Define $C_c(v)=\set{C_c^{s_i}(v):s_i\in S}$.
The set $C_c(v)$ can also be defined recursively by induction on the length of $v$ and the rank of $W$.
Recall now that $c$-sortable elements are the bottom elements of $c$-Cambrian congruence classes, and that each $c$-Cambrian congruence class defines a cone.
Surprisingly, the cone for $v$ is completely described by the data $C_c(v)$.
Specifically, by \cite[Theorem~6.3]{typefree}, the cone has $n$ facets, and the vectors in $C_c(v)$ are the normal vectors to these facets, pointing into the interior of the cone.

In the example of $B_2$, we refer to the Coxeter fan and $s_1s_2$-Cambrian fan pictured in Figure~\ref{B2Fc fig}.
The region $D$ is, as usual, identified with the identity element $1$.
Let $\alpha_1$ and $\alpha_2$ be the inward facing normal vectors to $D$, with $\alpha_1$ orthogonal to the line separating $D$ from the region identified with $s_1$.
Then $s_1$ is the reflection orthogonal to $\alpha_1$ and $s_2$ is the reflection orthogonal to $\alpha_2$.
Let us calculate $C_{s_1s_2}(s_2)$.
The letter $s_1$ is skipped in the first step of constructing the $(s_1,s_2)$-sorting word for $s_2$, so $C_{s_1s_2}^{s_1}(s_2)$ is $\alpha_1$.
In the next step, $s_2$ becomes the first letter of the $(s_1,s_2)$-sorting word.
The next time $s_2$ is tried, the sorting word is already complete, so $s_2$ is skipped.
Thus $C_{s_1s_2}^{s_2}(s_2)$ is $s_2\alpha_2=-\alpha_2$.
We see in Figure~\ref{B2Fc fig}.b. that $\alpha_1$ and $-\alpha_2$ are indeed inward-facing normal vectors to the cone associated to $s_2$.
Similarly, we calculate $C_{s_1s_2}(s_1s_2)$.
The letter $s_1$ is skipped after both $s_1$ and $s_2$ have already been tried and placed in the sorting word, so $C^{s_1}_{s_1s_2}(s_1s_2)=s_1s_2\alpha_1$.
For the same reason, $C^{s_2}_{s_1s_2}(s_1s_2)=s_1s_2\alpha_2$, and we verify that the cone associated to $s_1s_2$ is indeed defined by inward-facing normal vectors $s_1s_2\alpha_1$ and $s_1s_2\alpha_2$.

\section{Applications}\label{app sec}
In this section, we briefly mention several mathematical applications of Cambrian lattices and sortable elements.

\subsection{Coxeter-Catalan combinatorics}
The combinatorics of objects counted by the Catalan numbers has a rich history.
Recently, several of these \emph{Catalan objects} have been shown to be special cases (the case $W=S_{n+1}$) of combinatorial constructions whose input is a finite Coxeter group $W$.
The more general objects are counted by a generalization of the Catalan number.
For more background, including references, we refer the reader to \cite{Armstrong,rsga}.

One of the Coxeter-Catalan constructions, the generalized associahedron, has already been mentioned.
Its vertices are indexed by a generalization of triangulations.
Another construction generalizes the classical noncrossing partitions.
Initially, the $W$-noncrossing partitions were noticed to be equinumerous with the vertices of the $W$-associahedron, for any $W$, but no uniform explanation was known.
Sortable elements were used in~\cite{sortable} to give the first uniform explanation, by uniformly defining a bijection from sortable elements to noncrossing partitions and a bijection to the vertices of the generalized associahedron.
(Although the bijections were uniformly \textbf{defined} in~\cite{sortable}, the proofs of some key lemmas were non-uniform.  
These lemmas were later proved uniformly in~\cite{typefree}.
Meanwhile, another uniform bijection between noncrossing partitions and vertices of the generalized associahedron was given in~\cite{ABMW}, using ideas from~\cite{BWlattice}.)

Recall the discussion at the end of Section~\ref{sort sec} of three types of combinatorial traction provided by sortable elements.
The bijection to vertices of the generalized associahedron is defined in terms of sorting words for sortable elements, and both bijections are proved using induction on length and rank.
The bijections also make the natural search tree structure on sortable elements available in the context of noncrossing partitions and in the context of the generalized associahedron.

\subsection{Sortable elements in quiver theory}
Sortable elements and Cambrian lattices have figured in the study of quiver representations and various topics.
Details are far beyond the scope of this article, but we refer the reader to the research papers \cite{Amiot,AIRT,Ing-Tho,KP} and to Thomas' chapter \cite{Thomas} in this volume.

\subsection{Combinatorial Hopf algebras}
Besides the Cambrian construction, there is another profitable generalization of the Tamari lattices in a different direction, which nevertheless has something in common with the Cambrian picture:  lattice congruences.
In \cite{LR,LRorder}, Loday and Ronco constructed a combinatorial Hopf algebra whose operations are closely tied to the combinatorics of the Tamari lattice.
A permutations-to-triangulations map (in an equivalent formulation using planar binary trees rather than triangulations) was also relevant to the Loday-Ronco Hopf algebra.
An infinite family of Hopf algebras, containing the Loday-Ronco Hopf algebra, is constructed in~\cite{con_app}, building on the observation that the permutations-to-triangulations map is a lattice homomorphism.
The construction in~\cite{con_app} yields Hopf algebras built on congruence classes of permutations, and suggests the project of instead realizing the Hopf algebras in terms of natural combinatorial objects.
The goal is a description of the product and coproduct in terms of natural operations on the combinatorial objects, and a description of the quotient lattice as a partial order on the combinatorial objects, with cover relations given by combinatorial moves analogous to diagonal flips on triangulations.
Besides the motivating examples of the Loday-Ronco Hopf algebra and the Hopf algebra of noncommutative symmetric functions~\cite{GKLLRT}, two Hopf algebras built on \emph{rectangulations} (tilings of a rectangle by rectangles) have already been studied~\cite{Giraudo,rectangle,generic}, as well as a Hopf algebra of \emph{sashes}~\cite{sash}, certain tilings counted by the Pell numbers.

\subsection{Cluster algebras}
Cluster algebras were defined by Fomin and Zelevinsky in~\cite{ca1} and have since become the subject of intense research, in part because of their unexpected connections to a variety of mathematical areas.
Sortable elements and Cambrian lattices/fans have found applications in the structural theory of cluster algebras.
These applications are surprising \textit{a priori}, but less surprising once the connection between Cambrian lattices and generalized associahedra is known, because generalized associahedra are the combinatorial structure underlying cluster algebras of finite type.
To give the reader an idea of these applications, we start with the gentlest possible introduction to cluster algebras, giving only the ideas behind the definition.
The lecture notes~\cite{rsga}, while still gentle, do contain some actual definitions.
Many other expository works are listed in the Cluster Algebras Portal, maintained by Sergey Fomin and easily found online.

The setting for a cluster algebra is a field of rational functions in $n$ variables.
We start with elements $x_1,\ldots x_n$ of the field and some ``combinatorial data.''
This pair \[\set{\mathrm{Combinatorial\ Data}},\,(x_1,\ldots x_n)\] is called the initial \emph{seed}.
The rational functions $x_1,\ldots x_n$ are the \emph{cluster variables}.
There is an operation called \emph{mutation}, which takes a seed and gives a new seed.
The combinatorial data tells how to do mutation.
Mutation can take place any of $n$ ``directions,'' and it does two things:
It switches out one cluster variable, replacing it with a new one, and it alters the combinatorial data.
The result is a new seed.
The mutation operation is involutive, in the sense that mutating in the $k\th$ direction, and then mutating the result in the $k\th$ direction returns the original seed.
\begin{figure}
\scalebox{0.9}{
\includegraphics{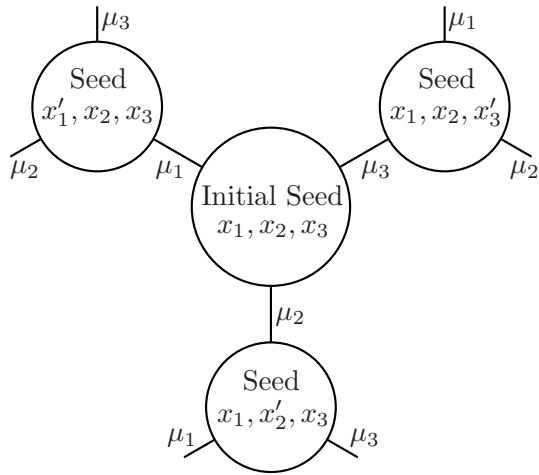}
\begin{picture}(0,0)
\put(-119,94){Initial Seed}
\put(-113.5,84){$x_1,x_2,x_3$}
\put(-104.5,29){Seed}
\put(-113.5,18){$x_1,x'_2,x_3$}
\put(-165,135){Seed}
\put(-175,124){$x'_1,x_2,x_3$}
\put(-43,135){Seed}
\put(-53,124){$x_1,x_2,x'_3$}
\put(-135,105){$\mu_1$}
\put(-93,53.5){$\mu_2$}
\put(-63,105){$\mu_3$}

\put(-186,105){$\mu_2$}
\put(-154,157){$\mu_3$}

\put(-32.5,157){$\mu_1$}
\put(-11,105){$\mu_2$}

\put(-131,11){$\mu_1$}
\put(-67,11){$\mu_3$}
\end{picture}
}
\caption{A schematic representation of seed mutations}
\label{exch fig}
\end{figure}
The basic setup is represented schematically in Figure~\ref{exch fig} for $n=3$, with dangling edges to indicate that the picture continues outward indefinitely.
Mutation in the $k\th$ direction is indicated by the symbol $\mu_k$.

To define the cluster algebra, we do all possible sequences of mutations, and collect all the cluster variables which appear.
The \emph{cluster algebra} for the given initial seed is the algebra of rational functions generated by all cluster variables.  
As a matter of bookkeeping, we place the seeds on the vertices of an infinite tree, with each vertex incident to $n$ edges labeled $\mu_1$ through $\mu_n$ to indicate the direction of mutation.
Once the initial seed is placed on some vertex of the tree, since each other vertex is connected to the initial vertex by a unique path, each vertex indicates a seed: the seed obtained from the initial seed by the sequence of mutations labeling the edges on this unique path.

However, depending on the initial combinatorial data, seeds at different vertices on the tree may coincide.
In this case, we obtain a smaller graph by identifying vertices of the infinite tree whenever they have the same seed.
The smaller graph is called the \emph{exchange graph}.
The cluster algebra is said to be of \emph{finite type} when the exchange graph is a finite graph.
For example, when $n=2$, the infinite tree is an infinite path.
There are four cases where the initial combinatorial data causes this infinite path to collapse to a cycle (of length $4$, $5$, $6$ or $8$).
For larger $n$, the main result of~\cite{ca2} identifies the exchange graphs of cluster algebras of finite type as the $1$-skeleta of generalized associahedra.

To see how Cambrian lattices and sortable elements contribute, it is first useful to know a bit more about the combinatorial data that describes the initial seed.
This data takes the form of an $n\times n$ matrix $B$, called the \emph{exchange matrix}, and some other data that should be thought of as a choice of coefficients.
The key point is that choosing $B$ is equivalent to choosing a Coxeter group $W$, a specific representation of $W$ as a group generated by reflections, and an orientation of the Coxeter diagram of $W$.
(Strictly speaking, we must choose a \emph{crystallographic} Coxeter group, thus allowing the representation of $W$ to be given by integer matrices.)
When $W$ is finite, the orientation of its diagram defines a Cambrian lattice, and having a specific representation of $W$ then implies a specific realization of the corresponding Cambrian fan.
The Hasse diagram of this Cambrian lattice is isomorphic to the exchange graph of the cluster algebra defined by $B$, with the initial seed corresponding to the bottom element of the Cambrian lattice (the identity element in $W$).
But in fact, even more structural information is contained in the sortable elements and Cambrian fan.
The cluster variables appearing in the various seeds can be identified by various statistics known as \emph{denominator vectors} and \emph{$\mathbf{g}$-vectors}.
The denominator vectors are easily read off from the sorting words of sortable elements.
The $\mathbf{g}$-vectors are encoded geometrically in the Cambrian fan, which coincides with what one might call the \emph{$\mathbf{g}$-vector fan} of the cluster algebra.
Thus we have arrived, via the lattice theory and combinatorics of Coxeter groups, at the fundamental combinatorial and polyhedral structure underlying cluster algebras.  
In particular, Figures~\ref{B2Fc fig}.b and~\ref{A3camb fig} are illustrations of this underlying structure.

\begin{figure}
\rotatebox{90}{\scalebox{.6}{\includegraphics{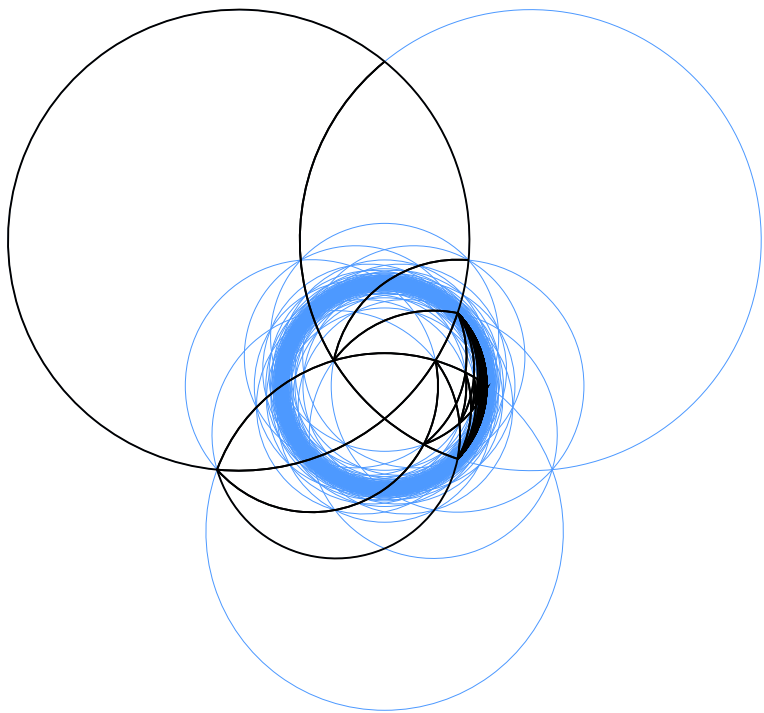}}}
\caption{A Cambrian fan for an infinite Coxeter group}
\label{infcamb fig}
\end{figure}

The definition of the weak order and the definition of sortable elements are both valid even when $W$ is not of finite type.
The weak order becomes a semilattice (i.e.\ meets exist) but not a lattice, so the lattice-theoretic definition of the Cambrian lattice is no longer available. 
However, the weak order restricted to sortable elements \textbf{is} still available and useful.
This \emph{Cambrian semilattice} serves as a model for \textbf{part of} the exchange graph, and in that part, we can still read off $\mathbf{g}$-vectors, exchange matrices and principal coefficients.
We can also define a Cambrian fan using the sets $C_c(v)$, and this Cambrian fan coincides with part of the $\mathbf{g}$-vector fan.
An example of an infinite Cambrian fan is shown in Figure~\ref{infcamb fig}, drawn in the same stereographic projection as previous fan pictures.
The cones of the Cambrian fan are outlined in black, and the reflecting hyperplanes are drawn as thin blue (or light gray) curves.
As one moves right from the center of the figure, one encounters infinitely many smaller and smaller Cambrian cones.

In an important infinite case, the Cambrian model still recovers the entire $\mathbf{g}$-vector fan.
When $W$ is of \emph{affine} type, the $\mathbf{g}$-vector fan is the union of a Cambrian fan with the antipodal opposite of a Cambrian fan for the opposite orientation of the diagram.
\begin{figure}
\rotatebox{90}{\scalebox{.6}{\includegraphics{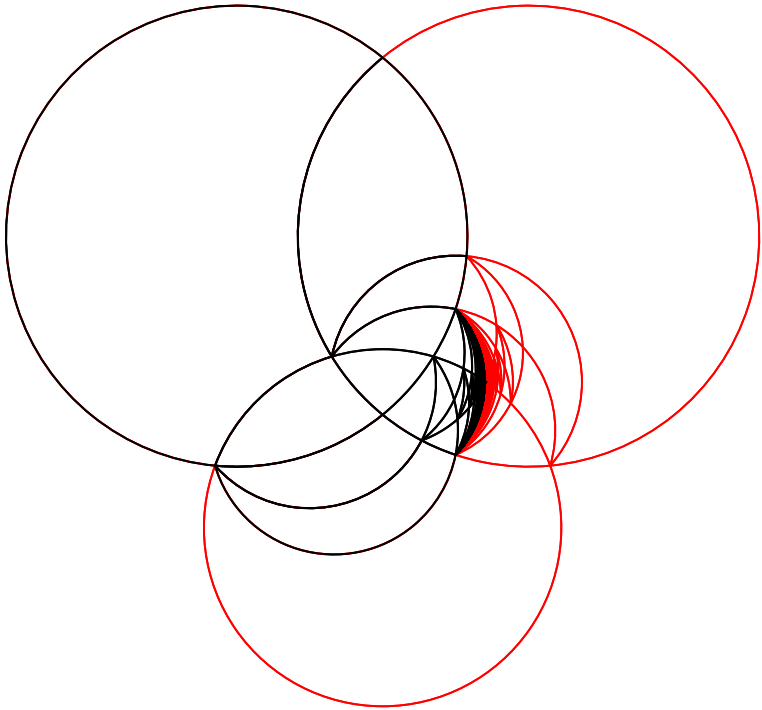}}}
\caption{A doubled Cambrian fan for an affine Coxeter group}
\label{infcambboth fig}
\end{figure}
Such a \emph{doubled Cambrian fan}, built on the fan from Figure~\ref{infcamb fig}, is shown in Figure~\ref{infcambboth fig}.
The second Cambrian fan is indicated by the red (or gray) arcs, some of which are obscured by the original black arcs.

The Cambrian model not only recovers the combinatorics and polyhedral geometry of cluster algebras, but also leads to proofs of some key structural conjectures about cluster algebras made in \cite{FZ-CDM,ca4}, in the cases where the underlying Coxeter group is of finite or affine type.
For some of these cases, the Cambrian model is currently the only known approach to proving the structural results.

Of necessity, we have given almost no details on cluster algebras and how they interact with Cambrian lattices and sortable elements.
Readers interested in filling in the details can find the necessary cluster algebras definitions in \cite{ca4}, details about $\mathbf{g}$-vectors and Cambrian fans in \cite{camb_fan,framework,YZ}, and more about the Cambrian model of cluster combinatorics in \cite{framework,aff_camb}.


\end{document}